%%%%%%%%%%%%%    Algebraic and Geometric Topology: agt-1-4.tex  %%%%%%%%
%%%%        
%%%%                   On asymptotic dimension of groups
%%%%             
%%%%                      G. Bell and A. Dranishnikov      
%%%%  
%%%%                 Published in Volume 1(2001) pages 57-71
%%%%
%%%%                   Publication date 27 January 2001
%%%%
%%%%                       This is a plain TeX file
%%%%
%%%%
%%%%%%%%%%%%%%%%%%                                   %%%%%%%%%%%%%%%%%%%
%%%%%%%%%%%%%%%%%%%%%%%%%%%%%%%%%%%%%%%%%%%%%%%%%%%%%%%%%%%%%%%
%%%%%%%%%%%             gtmacros.tex            %%%%%%%%%%%%%%%
%%%%%%%%%%%             version 1.6             %%%%%%%%%%%%%%% 
%
%                       Colin Rourke   
%
%
%    These macros are recommended for use by authors submitting articles   
%    to Geometry and Topology or to Algebraic and Geometric Topology.  
%    They are intended to be used with plain TeX. Each macro is described 
%    briefly to make it clear how to use it (or to modify it to achieve
%    different results).  If you modify this file then please change its
%    name.  If you modify this file and use the modified file to 
%    format an article for submission to Geometry and Topology or
%    Algebraic and Geometric Topology, then please paste the modified
%    file into your main TeX file.  Do not submit it as a separate file.
%      
%    Instructions on using these macros are also given in  gtmacins.tex  
%    or  gtmacins.ps  or .pdf  available on the gt www pages or by 
%    anonymous ftp from the gt/info/macros directory.
%
%
\magnification=\magstephalf      % Sets default point size to 11pt.
%
%  Basic layout parameters :
%
\vsize=7.5truein                 % Sets text height to 7.5 inches.
\hsize=5.2truein                 % Sets text width to 5.2 inches.
\newskip\stdskip                 % standard vertical space
\stdskip=6pt plus3pt minus3pt    % (slightly more stretchy
\medskipamount=\stdskip          % than the usual \medskip)
\parindent=0pt                   % Paragraphs are non-indented with
\parskip=\stdskip                % a little space between paragraphs. 
\abovedisplayskip=\stdskip       %  Reduces the space
\belowdisplayskip=\stdskip       %  around displays.
\mathsurround=0.75pt             % Gives a little extra space around maths.
\overfullrule=0pt                %  Prevents black boxes
%
%   The following macro is for principal paragraph breaks ie
%   a paragraph break with a slightly larger space :
%
\def\ppar{\par\goodbreak\vskip 8pt plus 4pt minus 4pt}     
%
%  The standard horizontal space for theorems, labels etc :
%
\def\stdspace{\hskip 0.75em plus 0.15em\ignorespaces}
\let\qua\stdspace % useful abbreviation (3/4 of a quad)
%
%%%%%%%%%%%%%%            FONT MACROS            %%%%%%%%%%%%
%
%           The following font macros define the AMS symbol 
%           and Euler-Fraktal fonts for use in text and
%           mathematics with appropriate size changes.
%           They also define two new control sequences  
%           \small  and  \large  (similar to those built
%           into LaTeX) which change the size of all fonts 
%           both in text and maths.  \small  is 10% smaller 
%           than normal and  \large  30% bigger.  The strange
%           size of the \large text fonts (10pt scaled 1315)
%           is because these macros are intended to be used
%           at \magstephalf.  The result is 10pt scaled 1440
%           (\magstep2) which is a standard font size.  If
%           you are borrowing these macros to use them at
%           another basic  \magnification, then you will
%           probably need to change 1315 to 1200 in the eleven
%           places marked ** below.  \large  will then be
%           20% bigger than normal.  Note that at \magstephalf
%           all the fonts come out roughly one point larger
%           than their size as defined in these macros.
%
%           The size-changing macros are based on Knuth's
%           \ninepoint and \eightpoint macros.
%
%
%    The macros are laid out in a way which makes it clear how to
%    add futher fonts (or delete unavailable ones) and how to add
%    further size changes.
%
%    First comes a definition of  \hexnumber  which is needed to
%    refer to font families whose family number is not known :
%
\def\hexnumber#1{\ifcase#1 0\or 1\or 2\or 3\or 4\or 5\or 6\or 7\or 8\or
 9\or A\or B\or C\or D\or E\or F\fi}
%
%     Next we define the AMS symbol-a fonts at 13,10,9,7,6,5 points
%
\font\thirtnmsa=msam10 scaled 1315    %%% **  see note above 
\font\tenmsa=msam10          \font\ninemsa=msam9
\font\sevenmsa=msam7         \font\sixmsa=msam6
\font\fivemsa=msam5
%%%%%%  (add further sizes here if you need them)
%
%    and the standard size family for these fonts
%
\newfam\msafam                  \textfont\msafam=\tenmsa
\scriptfont\msafam=\sevenmsa    \scriptscriptfont\msafam=\fivemsa
\edef\hexa{\hexnumber\msafam}        %  The msa family is  \fam\hexa
\def\msa{\fam\msafam\tenmsa}         %  \msa  switches to this family
%
%    Repeat these steps for the AMS symbol-b fonts
%
\font\thirtnmsb=msbm10 scaled 1315   %%%  ** see note above
\font\tenmsb=msbm10      \font\ninemsb=msbm9
\font\sevenmsb=msbm7     \font\sixmsb=msbm6
\font\fivemsb=msbm5
%%%%%%  (add further sizes here if you need them)
%
\newfam\msbfam                   \textfont\msbfam=\tenmsb       
\scriptfont\msbfam=\sevenmsb     \scriptscriptfont\msbfam=\fivemsb
\edef\hexb{\hexnumber\msbfam}    %  The msb family is \fam\hexb  
\def\msb{\fam\msbfam\tenmsb}     %  \msb switches to this family
%
%        Repeat for the Euler-Fraktal fonts 
%
\font\thirtneufm=eufm10 scaled 1315   %%% **  see note above 
\font\teneufm=eufm10                 \font\nineeufm=eufm9
\font\seveneufm=eufm7                \font\sixeufm=eufm6
\font\fiveeufm=eufm5
%%%%%%  (add further sizes here if you need them)
%
\newfam\eufmfam                    \textfont\eufmfam=\teneufm
\scriptfont\eufmfam=\seveneufm     \scriptscriptfont\eufmfam=\fiveeufm
\edef\hexf{\hexnumber\eufmfam}      % The Euler-Fraktal family is
\def\frak{\fam\eufmfam\teneufm}     % \fam\hexf and \frak switches to this
%
%%%  Add further fonts families here (using the same format) if you need
%    them.  The def of hexnumber is optional (it is only used for
%    \mathchardef 's).
%
%      Now we need to define the standard fonts (which are
%      already defined at 10,7 and 5 point) at 13,9 and 6 point:
%
%      Roman fonts:
\font\thirtnrm=cmr10 scaled 1315    %%%  ** see note above
\font\ninerm=cmr9                   \font\sixrm=cmr6   
%%%%%%  (add further sizes here if you need them)
%
%      Math italic fonts
\font\thirtni=cmmi10 scaled 1315    %%%  ** see note above 
\font\ninei=cmmi9                   \font\sixi=cmmi6  
%%%%%%  (add further sizes here if you need them)
%
%     Symbol fonts
\font\thirtnsy=cmsy10 scaled 1315   %%%  ** see note above
\font\ninesy=cmsy9                  \font\sixsy=cmsy6  
%%%%%%  (add further sizes here if you need them)
%
%     Bold face
\font\thirtnbf=cmbx10 scaled 1315   %%%  ** see note above 
\font\ninebf=cmbx9                  \font\sixbf=cmbx6  
%%%%%%  (add further sizes here if you need them)
%
%     The maths extension font (only defined at text size)
%
\font\thirtnex=cmex10 scaled 1315   %%%  ** see note above
\font\nineex=cmex9                  
%%%%%%  (add further sizes here if you need them)
%
%     Finally three fonts (text italic, slanted and typewriter type)
%     which are also only defined at text size
%
\font\thirtnit=cmti10 scaled 1315  %%%  ** see note above 
\font\nineit=cmti9                  
%%%%%%  (add further sizes here if you need them)
%
\font\thirtnsl=cmsl10 scaled 1315  %%%  ** see note above 
\font\ninesl=cmsl9                  
%%%%%%  (add further sizes here if you need them)
%
\font\thirtntt=cmtt10 scaled 1315  %%%  ** see note above 
\font\ninett=cmtt9                  
%%%%%%  (add further sizes here if you need them)
%
%
%     Now come the two main macros.  What  \small  does is to
%     change all the families of fonts from normal size which is
%     10,7,5  (ie 10pt text, 7pt subscript, 5pt subsubscript)
%     to 9,6,5.  \large  similarly changes to  13,9,7.  To make
%     other size changing macros, choose your three sizes, add
%     font size definitions if necessary and make the obvious changes
%     to one of these macros.  Change  \normalbaselineskip  and
%     \strutbox  dimensions to appropriate sizes as well.  To
%     add further fonts, insert them in each macro, using the
%     AMS fonts as a model.
%      
%
\def\small{%
%
%   redefine the sizes of the roman fonts :
%
\textfont0=\ninerm \scriptfont0=\sixrm \scriptscriptfont0=\fiverm
\def\rm{\fam0\ninerm}%       % ( \rm  sets \ninerm  in text mode
%                            %  and \fam0 in math mode)
%
%   and the math italic fonts :
%
\textfont1=\ninei \scriptfont1=\sixi \scriptscriptfont1=\fivei
%
%   and the symbol fonts :
%
\textfont2=\ninesy \scriptfont2=\sixsy \scriptscriptfont2=\fivesy
%
%   There is only one math extension font :
%
\textfont3=\nineex \scriptfont3=\nineex \scriptscriptfont3=\nineex
%
%   Next the bold font (named rather than numbered) :
%
\textfont\bffam=\ninebf \scriptfont\bffam=\sixbf
\scriptscriptfont\bffam=\fivebf \def\bf{\fam\bffam\ninebf}%
%
%   and the three text-only fonts : 
%
\textfont\itfam=\nineit \def\it{\fam\itfam\nineit}%
\textfont\slfam=\ninesl \def\sl{\fam\slfam\ninesl}%
\textfont\ttfam=\ninett \def\tt{\fam\ttfam\ninett}%
%
%   Now the three new families of AMS fonts :
%
%   AMS symbol-a
%
\textfont\msafam=\ninemsa \scriptfont\msafam=\sixmsa
\scriptscriptfont\msafam=\fivemsa \def\msa{\fam\msafam\ninemsa}%         
%
%   AMS symbol-b
%
\textfont\msbfam=\ninemsb \scriptfont\msbfam=\sixmsb
\scriptscriptfont\msbfam=\fivemsb \def\msb{\fam\msbfam\ninemsb}%         
%
%   Euler-Fraktal font
%
\textfont\eufmfam=\nineeufm  \scriptfont\eufmfam=\sixeufm
\scriptscriptfont\eufmfam=\fiveeufm \def\frak{\fam\eufmfam\nineeufm}%
%
%%%  Add further fonts families here if you need them.
%
%    Reset \normalbaselineskip and \strubox :
%
\normalbaselineskip=11pt%
\setbox\strutbox=\hbox{\vrule height8pt depth3pt width0pt}%
%
%    Set \normalbaselines and \rm (roman) as defaults :
%
\normalbaselines\rm
%
%    Reset some of the basic vertical skips:
%
\stdskip=4pt plus2pt minus2pt    
\medskipamount=\stdskip          
\parskip=\stdskip                
\abovedisplayskip=\stdskip       
\belowdisplayskip=\stdskip       
\def\ppar{\par\goodbreak\vskip 6pt plus 3pt minus 3pt}%     
%
%   And finally reset the size of section heads (see below):
%
\def\section##1{\global\advance\sectionnumber by 1
\vskip-\lastskip\penalty-800\vskip 20pt plus10pt minus5pt 
\egroup{\bf\number\sectionnumber\quad##1}\bgroup\small         
\vskip 6pt plus3pt minus3pt
\nobreak\resultnumber=1}%      % Reset resultnumber at start of section
}    %%%   End of  \small  macro      
%
%   Two useful abbreviations to keep track of \small material:
\def\beginsmall{\bgroup\small}
\let\endsmall\egroup
%
%
%    The \large  macro is similar (comments abbreviated):
%
%
\def\large{%
\textfont0=\thirtnrm \scriptfont0=\ninerm \scriptscriptfont0=\sevenrm
\def\rm{\fam0\thirtnrm}%
\textfont1=\thirtni \scriptfont1=\ninei \scriptscriptfont1=\seveni
\textfont2=\thirtnsy \scriptfont2=\ninesy \scriptscriptfont2=\sevensy
\textfont3=\thirtnex \scriptfont3=\thirtnex \scriptscriptfont3=\thirtnex
\textfont\bffam=\thirtnbf \scriptfont\bffam=\ninebf
\scriptscriptfont\bffam=\sevenbf \def\bf{\fam\bffam\thirtnbf}%
\textfont\itfam=\thirtnit \def\it{\fam\itfam\thirtnit}%
\textfont\slfam=\thirtnsl \def\sl{\fam\slfam\thirtnsl}%
\textfont\ttfam=\thirtntt \def\tt{\fam\ttfam\thirtntt}%
%   AMS symbol-a  :
\textfont\msafam=\thirtnmsa \scriptfont\msafam=\ninemsa
\scriptscriptfont\msafam=\sevenmsa \def\msa{\fam\msafam\thirtnmsa}%         
%   AMS symbol-b  :
\textfont\msbfam=\thirtnmsb \scriptfont\msbfam=\ninemsb
\scriptscriptfont\msbfam=\sevenmsb \def\msb{\fam\msbfam\thirtnmsb}%         
%   Euler-Fraktal font :
\textfont\eufmfam=\thirtneufm  \scriptfont\eufmfam=\nineeufm
\scriptscriptfont\eufmfam=\seveneufm \def\frak{\fam\eufmfam\teneufm}%
%%%% Add further fonts families here if you need them.
%   Reset \normalbaselineskip and \strubox and initialise :
\normalbaselineskip=16pt%
\setbox\strutbox=\hbox{\vrule height11.5pt depth4.5pt width0pt}%
\normalbaselines\rm}%
\let\Large\large   %  for compatibility with latex
%
%   The next two lines define commonly used switches for
%   blackboard bold (\Bbb) and gothic type (\goth).  The   
%   \Bbb  switch is set to work in the same way as in amstex
%   and switches only the next character to blackboard bold.

%
%   To use the new AMS fonts you can either use the control
%   sequences \msa, \msb (alias \Bbb) and \frak (alias \goth) eg :

   % see the msam font table
%
%   or, more generally, make \mathchardef's (cf Knuth p155) eg :
\mathchardef\plussquare="0\hexa01
\mathchardef\nge="3\hexb0B
\mathchardef\maltesecross="0\hexa7A
\mathchardef\del="0\hexf01
%
%   or you can use the amstex names for all the new symbols by
%   inserting the line  \input amsnames  in your file directly
%   after \input gtmacros. 
%   This presupposes that you have collected a copy of the file
%   amsnames.tex  from the  gt/info/macros  ftp directory.
%
%
%   Finally we need a small capital font (for author(s)) :
%
\font\sc=cmcsc10
%
%%%%%%%%%%%%%%%%%       END OF FONT MACROS     %%%%%%%%%%%%%
%
%
%                 Knuth's \square macro :
%
\def\sqr#1#2{{\vcenter{\vbox{\hrule  height.#2truept
	\hbox{\vrule width.#2truept height#1truept 
	\kern#1truept \vrule width.#2truept}
	\hrule height.#2truept}}}}
\def\sq{\sqr55}    %   A small square for end-of-proofs. 
%                  %   (Define other size squares by varing the
%                  %   the two numbers.)
%
%
%      Style macros for section heads, theorem statements etc :
%   
%
\newcount\sectionnumber            %%%  Allocate registers to take
\newcount\resultnumber             %%%  section and result numbers.
\sectionnumber=0\resultnumber=1    %%%  Set these registers to 0 and 1
%
%   The \section macro produces a \large bold faced section heading
%   numbered to the left.  Pagebreaks are encouraged before the
%   start of the section and discouraged directly after the heading.
%   Typical use  \section{First steps}  with typical result :
%
%    1  First Steps     (set bold and \large)
%
\def\section#1{\global\advance\sectionnumber by 1
\xdef\nextkey{\number\sectionnumber}%      (used by the \key macro)
\vskip-\lastskip\penalty-800\vskip 20pt plus10pt minus5pt 
{\large\bf\number\sectionnumber\quad#1}         
\vskip 8pt plus4pt minus4pt
\nobreak\resultnumber=1}      % Reset resultnumber at start of section
%
%
%
%   Next a macro to set subheadings (like the  \section  macro
%   but without the number, with less space and set in standard size).
%
%   Typical use :  \sh{Example formats}
%
         
%
%   The \proc ... \endproc macros ("proclaim") are for setting theorems, 
%   lemmas, conjectures etc with automatic numbering.  Typical use :    
%  
%    \proc{Theorem}Every lemon is yellow.\endproc
%
%   Typical result :
%     
%    Theorem 3.4  Every lemon is yellow.   

%   (with Theorem 3.4 set bold and a \stdspace of space before the 
%   statement set in slanted type).
%
\def\proc#1{\xdef\nextkey{\number\sectionnumber.\number\resultnumber}%
\vskip-\lastskip\ppar\bf%
\noindent#1\ \number\sectionnumber.\number\resultnumber
\stdspace\sl\global\advance\resultnumber by 1\ignorespaces}
 
%
%  The \prf ... \endprf macros are for setting proofs.  The code
%  for \prf includes the code for \endproc, so there is no need to
%  type \endproc if the theorem is followed immediatedly by a proof.
%
\def\prf{\vskip-\lastskip\ppar\noindent{\bf Proof}%
\stdspace\rm}                            %  For start of proofs  
   %  For end (or absence) of proofs
\def\endprf{\unskip\stdspace\hbox{}%     %  For end of proof (with
\hfill$\sq$\par\medskip}                 %  extra vertical space)  
        %  For start of proof with alternative name
              %  \endproof is an alias for \endprf
%
%   Typical uses :    
%  
%    \proc{Theorem}Every lemon is yellow. \qed\endproc
%
%    \proc{Theorem}Every lemon is yellow.
%    \prf Use your eyes. \endprf
%
%    \proc{Theorem}Every lemon is yellow.
%    \proof{Proof of theorem} Use your eyes. \endprf
%
%   The next macro is a variant of the \proc macro.  It has
%   exactly the same result except that it omits the number.
%
%   Typical use :  
%    
%    \proclaim{Conjecture}Some oranges are yellow.\endproc
%
\def\proclaim#1{\vskip-\lastskip\ppar\bf%
\noindent#1\stdspace\sl\ignorespaces} 
\let\endproclaim\endproc
%
%   The next macro is a further variant for remarks, definitions etc.   
%   It omits the number and does not switch on slanted type.  
%  
%   Typical use :
%
%    \rk{Remark}Some lemons are thick-skinned.\endrk
%
\def\rk#1{\vskip-\lastskip\ppar{\bf #1}\stdspace\ignorespaces}                

%
%   The next macro is for numbering equations etc, \label  produces the 
%   correct number  x.y  and advances the resultnumber register
%
%   Typical use :
%
%     $$fx=7\eqno{\bf\label}$$
%
%   result :
%
%                           fx = 7                           3.5
%
\def\label{\xdef\nextkey{\number\sectionnumber.\number\resultnumber}%
\number\sectionnumber.\number\resultnumber
\global\advance\resultnumber by 1}
%
%
%
%   The next macros are to automate external references.  To use them 
%   type \reflist ..... \endreflist near the beginning of your paper, 
%   where  .... is the list of references in alphabetical order 
%   and in  the form  \key{KEY}  reference    where "KEY" is a 
%   string of characters which reminds you of the reference.   
%   Separate  references with a blank line or a \par.   Eg 
%
%     \reflist
%
%     ..... more references ....
%
%     \key{Kn-84} {\bf D Knuth}, {\it The TeXbook}, Addison--Wesley (1984)
%
%     ..... more references ....
%
%     \endreflist
%
%   Then type  \references  where you wish the references to be printed
%   (normally near the end of the paper).  To refer to Knuth type
%   for example    see Knuth [\ref{Kn-84}, page 320]   and the correct
%   numerical reference will be printed.  Edit the \references macro
%   to change the formatting (if desired).
%   There is an alternative \refkey for \key, provided your KEY contains
%   only letters.  The syntax is:
%
%     \reflist
%
%     ..... more references ....
%
%     \refkey\Knuth  {\bf D Knuth}, {\it The TeXbook}, Addison--Wesley (1984)
%
%     ..... more references ....
%
%     \endreflist
%
%   \key{Knuth}  has exactly the same maening as \refkey\Knuth and you
%   can mix the two syntaxes if you want.  But \refkey\Kn-84
%   would not work.  It would set Kn as the KEY and -84 would get printed!
%
\newcount\refnumber              %  Register for reference numbers
\refnumber=1                     %  set initially to 1.
\long\def\reflist#1\endreflist{%
\long\def\thereflist{#1}{\def\refkey##1##2\par{\xdef##1{\number\refnumber}%
\global\advance\refnumber by 1}%
\def\key##1##2\par{\expandafter\xdef%
\csname##1\endcsname{\number\refnumber}%
\global\advance\refnumber by 1}#1\par}}
\long\def\references{%
\penalty-800\vskip-\lastskip\vskip 15pt plus10pt minus5pt 
{\large\bf References}\ppar %`References' is set \large bold with space around.
{\leftskip=25pt\frenchspacing    % The list of references is set 
\small\parskip=3pt plus2pt       % \small  with small spaces between,
\def\refkey##1##2\par{\noindent  % numbers in [,]'s and set just to the
\llap{[##1]\stdspace}\ignorespaces##2\par}         % left of a 25pt margin.
\def\key##1##2\par{\noindent  
\llap{[\ref{##1}]\stdspace}\ignorespaces##2\par}  
\def\,{\thinspace}\thereflist\par}}
%
%   Next a footnote macro (with automatic numbering) which sets the
%   footnote  \small.
%
%   Typical use :
%         ..... are yellow.\fnote{By yellow here we mean Britsh
%    Standard colour BS3320.} 
%
\newcount\footnotenumber         % Register for footnote number
\footnotenumber=1                % set initially to 1
\def\fnote#1{\xdef\nextkey{\number\footnotenumber}%
{\small\ifnum\footnotenumber>9\parindent=14pt%
\else\parindent=10pt\fi\footnote{$^{\number\footnotenumber}$}%
{\hglue-5pt#1}\global\advance\footnotenumber by 1}}
%
%
%   Next macros for handling figures with automatic numbering (using 
%   TeX's \midinsert to float the figure to a suitable place).
%   
%   The \figure ... \endfigure macro centres the figure and adds
%   an automatically numbered label  Figure XX  after it.
%
%   If you have a caption, then type \caption{caption text} 
%   somewhere between \figure and \endfigure.  The macro
%   will then add  Figure XX: caption text  after the figure.
%
%   If you want an unnumbered or uncentred figure, then use TeX's raw 
%       \midinsert Figure instructions \endinsert  
%   and if you want a numbered figure label in the same style then
%   use \caption{caption text} outside of  \figure ... \endfigure.
%
%   If you need just the label Figure XX  outside of  \figure ... \endfigure
%   then type  \figurelabel .
%
\newcount\figurenumber          % register for figure number
\figurenumber=1                 % set initially to 1
\def\caption#1{\xdef\nextkey{\number\figurenumber}%
\cl{\small Figure \number\figurenumber: #1}%
\global\advance\figurenumber by 1}
\def\figurelabel{\xdef\nextkey{\number\figurenumber}%
\cl{\small Figure \number\figurenumber}%
\global\advance\figurenumber by 1}
\long\def\figure#1\endfigure{{\xdef\nextkey{\number\figurenumber}%
\let\captiontext\relax\def\caption##1{\xdef\captiontext{##1}}%
\midinsert\cl{\ignorespaces#1\unskip\unskip\unskip\unskip}\vglue6pt\cl{\small 
Figure \number\figurenumber\ifx\captiontext\relax\else: \captiontext
\fi}\endinsert\global\advance\figurenumber by 1}}
%
%   Macros for self-correcting internal references.
%
%   There are two macros  \key{KEY}  and  \ref{KEY} .
%
%   The \key macro sets up KEY as a key for whatever number is 
%   being referenced and the \ref macro converts the KEY into 
%   that number.  Type \key after a  \section or \proc or 
%   \label or \fnote or \figure or \caption or \figurelabel .
%
%   Example:
%
%       \section{Introduction}\key{intro}
%       \proc{Theorem}\key{MainTh}Lemons are yelloy\endproc
%       Here we follow\fnote{Follow in the sense of Dickens}
%       \key{Dickens-note}the crowd ....  
%
%       In section \ref{intro}
%       we stated theorem \ref{mainTh} and noted (see footnote 
%       \ref{Dickens-note}) ...
%
\def\nextkey{??}   %  initialise \nextkey (which is reset by all the
%                     numbering macros)
%
\def\key#1{\expandafter\xdef\csname #1\endcsname{\nextkey}}
\def\ref#1{\expandafter\ifx\csname #1\endcsname\relax
\immediate\write16{Reference {#1} undefined}??\else
\csname #1\endcsname\fi}
%
%   Note:  If the KEY contains only letters then \KEY has exactly the
%   same meaning as \ref{KEY} so in the example you could have:
%
%       In section \intro\ we ....
%
%   The \key will work at any time after the macro which sets the
%   number, provided no other macro which sets a number has been used. 
%
%   Macros for forward references:
%              =======
%   The \key \ref macros ONLY work for backwards references.  If you  
%   want to use forwards references, then type \useforwardrefs  near
%   the beginning of your file.  The KEY's are then stored in an
%   auxiliary  .ref  file and you then suffer the same disadvantage as
%   when using LaTeX that you must TeX the file twice to get
%   the references correct.
%
%   To use a forward ref type \ref{KEY}.  (You can type the
%   alternative  \KEY  but you'll get an error on first TeX'ing 
%   if the \KEY is not yet defined.) 
%
%   The macro also allows external references to be listed at the end 
%   of the file (if you wish to).  (Indeed they can be typed anywhere
%   before the \references command.)  You can combine the reference list
%   and the \references command by typing the references (using the
%   same syntax as before) between the commands \biblio and \endbiblio 
%   (don't type \references or they'll be printed twice).
%
\newread\gtinfile
\newwrite\gtreffile
\def\useforwardrefs{
\openin\gtinfile\jobname.ref
\ifeof\gtinfile
\closein\gtinfile
\immediate\write16{No file \jobname.ref}
\else
\closein\gtinfile
\input \jobname.ref
\fi
\immediate\openout\gtreffile \jobname.ref
%
%   Adapt \key :
%
\def\key##1{{\def\\{\noexpand}%
\expandafter\xdef\csname ##1\endcsname{\nextkey}%
\immediate\write\gtreffile{\\\expandafter\\\def\\\csname ##1\\\endcsname%
{\nextkey}}}}
%
%  Adapt macros for external references:  
%
\long\def\reflist##1\endreflist{%
\long\def\thereflist{##1}{\def\refkey####1####2\par{\xdef####1{%
\number\refnumber}{\def\\{\noexpand}\immediate\write\gtreffile
{\\\def\\####1{\number\refnumber}}}\global\advance\refnumber by 1}%
\def\key####1####2\par{\expandafter\xdef%
\csname####1\endcsname{\number\refnumber}%
{\def\\{\noexpand}\immediate\write\gtreffile
{\\\expandafter\\\def\\\csname ####1\\\endcsname{\number\refnumber}}}
\global\advance\refnumber by 1}##1\par}}
\long\def\biblio##1\endbiblio{\reflist##1\endreflist\references}%
%
%  Adapt obselete key macros (\numkey, \seckey and \figkey):
%
\def\numkey##1{{\def\\{\noexpand}%
\xdef##1{\number\sectionnumber.\number\resultnumber}
\immediate\write\gtreffile{\\\def\\##1%
{\number\sectionnumber.\number\resultnumber}}}}
\def\seckey##1{{\def\\{\noexpand}\xdef##1{\number\sectionnumber}
\immediate\write\gtreffile{\\\def\\##1{\number\sectionnumber}}}}
\def\figkey##1{\xdef##1{\number\figurenumber}%
{\def\\{\noexpand}\immediate\write\gtreffile%
{\\\def\\##1{\number\figurenumber}}}
\number\figurenumber\global\advance\figurenumber by 1}
}   %  end of \useforwardrefs
%
%
%   The next five macros are obselete and have been superseeded by
%   the general \key macro above.  They are included merely to 
%   maintain backward compatibility for the package:
%
%
\def\figkey#1{\xdef#1{\number\figurenumber}%
\number\figurenumber\global\advance\figurenumber by 1}
\def\fig#1#2\endfig{%
\midinsert\cl{#2}\vglue6pt\cl{\small Figure #1}\endinsert}
\def\newfig{\number\figurenumber\global\advance\figurenumber by 1}
\def\numkey#1{\xdef#1{\number\sectionnumber.\number\resultnumber}}
\def\seckey#1{\xdef#1{\number\sectionnumber}}
%
%   End of obselete macros.
%
%
%   The next macro is a version of the verbatim macro given by Knuth.
%
%   This macro produces a "verbatim" printout of
%   any ASCII string which does not contain the symbol "
%   (TeX files do not usually contain " 's).
%   More precisely, everything between consecutive pairs
%   of " 's is printed verbatim in the typewriter font cmtt.
%   For an explanation of how the macro works, see Knuth pp 420-1.
%
%   There are two switches: \verb (which switches the macro on)
%   and \brev which switches the macro off (the default).  When
%   the macro is switched off the symbol " has its usual 
%   meaning for TeX.  To use the macro, type \verb before use
%   and the use " to switch verbatim on and off.  Be careful
%   not to use " for any other purpose.  There is no need to
%   switch the macro off again unless you need to use " for
%   some other purpose (eg making  \mathchardef 's).  Note 
%   that the macro MUST BE OFF before inputting  amsnames.tex .
%
%   Whether the macro is on or off you can always use the
%   control sequence \dq (double quote) for " e.g.
%   \mathchardef\sum=\dq1350  is perfectly valid.
%   The control sequence \ttq is an abbreviation for
%   {\tt\dq}.  Thus "\ttq" will produce " (in cmtt)
%   inside a verbatim quote.
%
%
   %  define a code for " so it can be used when \verb is on
  %  code for " in cmtt
%
\def\verb{\catcode`\"=\active}       %  The main
\def\brev{\catcode`\"=12}            %  switches.
\brev                                %  Prime switches and
\verb                                %  switch on.
{\obeyspaces\gdef {\ }}              
{\catcode`\`=\active\gdef`{\relax\lq}}
\def"{%
\begingroup\baselineskip=12pt\def\par{\leavevmode\endgraf}%
\tt\obeylines\obeyspaces\parskip=0pt\parindent=0pt%
\catcode`\$=12\catcode`\&=12\catcode`\^=12\catcode`\#=12%
\catcode`\_=12\catcode`\~=12%
\catcode`\{=12\catcode`\}=12\catcode`\%=12\catcode`\\=12%
\catcode`\`=\active\let"\endgroup}
\brev      %   Finally switch the macro off (for safety)
%
%   Macros for itemised lists.   Typical use :
%    
%    \items
%    \item{(i)}Colours must be defined.
%    \item{(ii)}Colour cards may not be cited.
%    \enditems
%
%   Result :
%
%    (i)  Colours must be defined. 
%   (ii)  Colour cards may not be cited.
%
%
           % Start of itemised list         
         % end of itemised list   
\def\item#1{\par\leavevmode\llap{#1\stdspace}%
\ignorespaces}                             % labelled item
               % bulleted item.
%
%   The \quote ... \endquote macros are for typesetting quotations :
%

%
%   A few useful abbreviations :
%
    %  Colon with correct spacing for maps.
\def\np{\vfil\eject}         %  Forced page break (new page).
\def\nl{\hfil\break}         %  New line.
\def\cl{\centerline}         %  Centerline
        %  The journal title in recommended style
    %  for monographs
\def\agt{{\mathsurround=0pt\it$\cal A\mskip-.7mu$lgebraic \&\ 
$\cal G\mskip-2mu$eometric $\cal T\!\!$opology}}  % AGT
%
%    Finally some macros for automatic title page or header generation.
%    To use them type your header information using the following  
%    example as a guide :
%
%    Note that \\ is used as standard separator (for lines in \title and
%    \address, between authors and between email addresses or URL's)
%    and that \email, \url and \secondaddress are optional.
%

% Example:  \title{A short spoof paper\\with a two-line title}
% =======   \authors{Albert Einstein\\Leonardo da Vinci}
%           \address{IAS\\Princeton}\secondaddress{Renaissance\\Venice}
%           \email{ae@ias.princeton.edu\\ldv@ren.ven.hist}
%           \abstract 
%           A short spoof paper with a very short abstract.
%           \endabstract 
%           \primaryclass{00-01, 00-02}\secondaryclass{68-00, 68-01}
%           \keywords{Short, spoof, paper}
%           \maketitlepage
%
%
%    The title page or header will then be generated automatically.
%
%
%    Define the various ingredients of the title page:
%
\def\title#1{\def\thetitle{#1}}

\def\author#1{\edef\previousauthors{\theauthors}
 \ifx\theauthors\relax\def\theauthors{#1}\else
 \def\theauthors{\previousauthors\par#1}\fi}

\let\authors\author        % aliases
\def\address#1{\edef\previousaddresses{\theaddress}
 \ifx\theaddress\relax\def\theaddress{#1}\else
 \def\theaddress{\previousaddresses\par\vskip 2pt\par#1}\fi}
                             % alias
\def\secondaddress#1{\edef\previousaddresses{\theaddress}
 \ifx\theaddress\relax\def\theaddress{#1}\else
 \def\theaddress{\previousaddresses\par{\rm and}\par#1}\fi}   

\def\email#1{\edef\previousemails{\theemail}
 \ifx\theemail\relax\def\theemail{#1}\else
 \def\theemail{\previousemails\hskip 0.75em\relax#1}\fi}
  % aliases
\def\secondemail#1{\edef\previousemails{\theemail}
 \ifx\theemail\relax\def\theemail{#1}\else
 \def\theemail{\previousemails\hskip 0.75em{\rm and}\hskip 0.75em
 \relax#1}\fi}
\def\url#1{\edef\previousurls{\theurl}
 \ifx\theurl\relax\def\theurl{#1}\else
 \def\theurl{\previousurls\hskip 0.75em\relax#1}\fi}
      % aliases
\def\secondurl#1{\edef\previousurls{\theurl}
 \ifx\theurl\relax\def\theurl{#1}\else
 \def\theurl{\previousurls\hskip 0.75em{\rm and}\hskip 0.75em
 \relax#1}\fi}
\long\def\abstract#1\endabstract{\long\def\theabstract{#1}}
\def\primaryclass#1{\def\theprimaryclass{#1}}
                        % alias
\def\secondaryclass#1{\def\thesecondaryclass{#1}}
\def\keywords#1{\def\thekeywords{#1}}
%
%  Set \\ to \par and title page items to \relax to initialise macros :
%
\let\\\par\let\thetitle\relax\let\theshorttitle\relax
\let\theauthors\relax\let\theshortauthors\relax
\let\theaddress\relax\let\theshortaddress\relax
\let\theemail\relax\let\theurl\relax
\let\theabstract\relax\let\theprimaryclass\relax
\let\thesecondaryclass\relax\let\thekeywords\relax
%
%
%
%   Basic title page layout (edit this macro if you
%   wish to adjust the title page layout) :
%
\long\def\maketitlepage{    % start of definition of \maketitlepage

\vglue 0.2truein   % top margin

% title :
%
{\parskip=0pt\leftskip 0pt plus 1fil\def\\{\par\smallskip}{\large
\bf\thetitle}\par\medskip}   

\vglue 0.15truein 

% authors :
%
{\parskip=0pt\leftskip 0pt plus 1fil\def\\{\par}{\sc\theauthors}
\par\medskip}%
 
\vglue 0.1truein 

% address(es) email's and URL's (with switches to detect whether the
% optional items have been used) :
%
{\small\parskip=0pt
{\leftskip 0pt plus 1fil\def\\{\par}{\sl\theaddress}\par}
\ifx\theemail\relax\else  % email address?
\vglue 5pt \def\\{\stdspace{\rm and}\stdspace} 
\cl{Email:\stdspace\tt\theemail}\fi
\ifx\theurl\relax\else    % URL given?
\vglue 5pt \def\\{\stdspace{\rm and}\stdspace} 
\cl{URL:\stdspace\tt\theurl}\fi\par}

\vglue 7pt 

{\bf Abstract}

\vglue 5pt

\theabstract

\vglue 7pt 

{\bf AMS Classification numbers}\quad Primary:\quad \theprimaryclass\par

Secondary:\quad \thesecondaryclass

\vglue 5pt 

{\bf Keywords:}\quad \thekeywords

\np  % page break at the end of the title page

}    % end of definition of \maketitlepage
%
%    % \makeshorttitle (for general preprints) doesn't take a new page
%
\long\def\makeshorttitle{    % start of definition of \makeshorttitle

%\vglue 0.2truein   % top margin

% title :
%
{\parskip=0pt\leftskip 0pt plus 1fil\def\\{\par\smallskip}{\large
\bf\thetitle}\par\medskip}   

\vglue 0.05truein 

% authors :
%
{\parskip=0pt\leftskip 0pt plus 1fil\def\\{\par}{\sc\theauthors}
\par\medskip}%
 
\vglue 0.03truein 

% address(es) email's and URL's (with switches to detect whether the
% optional items have been used) :
%
{\small\parskip=0pt
{\leftskip 0pt plus 1fil\def\\{\par}{\sl\ifx\theshortaddress\relax
\theaddress\else\theshortaddress\fi}\par}
\ifx\theemail\relax\else  % email address?
\vglue 5pt \def\\{\stdspace{\rm and}\stdspace} 
\cl{Email:\stdspace\tt\theemail}\fi
\ifx\theurl\relax\else    % URL given?
\vglue 5pt \def\\{\stdspace{\rm and}\stdspace} 
\cl{URL:\stdspace\tt\theurl}\fi\par}

\vglue 10pt 

% abstract and classification numbers (with switches):

{\small\leftskip 25pt\rightskip 25pt{\bf Abstract}\stdspace\theabstract

{\bf AMS Classification}\stdspace\theprimaryclass
\ifx\thesecondaryclass\relax\else; \thesecondaryclass\fi\par
{\bf Keywords}\stdspace \thekeywords\par}
\vglue 7pt
}    % end of definition of \makeshorttitle
\let\maketitle\makeshorttitle        %% alias
%
%    %%%% \makeagttitle (for AGT) similar to \makeshorttitle but
%         with addresses omitted (they go at the end)
%
%%%% publication info and test defaults:

\def\volumenumber#1{\def\thevolumenumber{#1}}
\def\volumeyear#1{\def\thevolumeyear{#1}}
\def\pagenumbers#1#2{\def\startpage{#1}\def\finishpage{#2}}
\def\published#1{\def\publishdate{#1}}
\def\received#1{\def\receiveddate{#1}}
\def\revised#1{\def\reviseddate{#1}}
\let\reviseddate\relax
%% Defaults for authors to use to check layout
\volumenumber{X}
\volumeyear{20XX}
\pagenumbers{1}{XXX}
\published{XX Xxxember 20XX}

\long\def\makeagttitle{   %%% start of definition of \makeagttitle
\agt\hfill      %   Journal title (top left) 
%   logo placeholder (top right)
\hbox to 60truept{\vbox to 0pt{\vglue -14truept{\bf [Logo here]}\vss}\hss}
\break
{\small Volume \thevolumenumber\ (\thevolumeyear)
\startpage--\finishpage\nl
Published: \publishdate}

\vglue .2truein

% title
{\parskip=0pt\leftskip 0pt plus 1fil\def\\{\par\smallskip}{\large
\bf\thetitle}\par\medskip}   
\vglue 0.05truein 

% authors :
%
{\parskip=0pt\leftskip 0pt plus 1fil\def\\{\par}{\sc\theauthors}
\par\medskip}%
 
\vglue 0.03truein 

%  abstract and classification numbers:

{\small\leftskip 25truept\rightskip 25truept{\bf Abstract}\stdspace\theabstract

{\bf AMS Classification}\stdspace\theprimaryclass
\ifx\thesecondaryclass\relax\else; \thesecondaryclass\fi\par
{\bf Keywords}\stdspace \thekeywords\par}\vglue 7truept

}   %%%% end of definition of \makeagttitle

%%%%% Macro to typeset addresses (typically at the end of the paper)

\def\Addresses{\bigskip
{\small \parskip 0pt \leftskip 0pt \rightskip 0pt plus 1fil \def\\{\par}
\sl\theaddress\par\medskip \rm Email:\stdspace\tt\theemail\par
\ifx\theurl\relax\else\smallskip \rm URL:\stdspace\tt\theurl\par\fi}}

\def\agtart{%   Full mock-up of AGT article style (for authors to test with)
%  get print centerpage:
\hoffset 14truemm
\voffset 31truemm
\font\phead=cmsl9 scaled 950
\font\pnum=cmbx10 scaled 913
\font\pfoot=cmsl9 scaled 950
%  headline and footline
\headline{\vbox to 0pt{\vskip -4.5mm\line{\small\phead\ifnum
\count0=\startpage ISSN numbers are printed here
\hfill {\pnum\folio}\else\ifodd\count0\def\\{ }% 
\ifx\theshorttitle\relax\thetitle\else\theshorttitle\fi\hfill{\pnum\folio}
\else\def\\{ and }{\pnum\folio}\hfill\ifx\theshortauthors\relax\theauthors
\else\theshortauthors\fi\fi\fi}\vss}}
\footline{\vbox to 0pt{\vglue 0mm\line{\small\pfoot\ifnum\count0=\startpage
Copyright declaration is printed here\hfill\else
\agt, Volume \thevolumenumber\ (\thevolumeyear)\hfill\fi}\vss}}
%  force \agttitle
\let\maketitle\makeagttitle\let\makeshorttitle\makeagttitle}

%%%
%%%  This is agtout.tex.  
%%%
%%%  This the version of  gtoutput.tex  intended to finish formatting
%%%  papers published in Algebriac & Geometric Topology and stored in the
%%%  arXiv.   All versions of  gtoutput.tex  are copyright 
%%%  GT Publications and are to be used _only_ for formatting
%%%  the officially published version of ABT or G&T papers.
%%%
%%%
%%%                                             Colin Rourke  27.102000
%%%
%%%  To create header file  head.xxx  comment out the first \endinput

%  test for latex or plain tex
\def\ifplaintex{\expandafter\ifx\csname documentclass\endcsname\relax}

\def\gtp{{\mathsurround=0pt\it $\cal G\mskip-2mu$eometry \&\ 
$\cal T\!\!$opology $\cal P\!$ublications}}  % GT publications

\def\recd{{\small Received:\qua\receiveddate\ifx\reviseddate\relax
\else\qquad Revised:\qua\reviseddate\fi\par}} 

%  define the various new ingredients of the title page and the data
%  output files

\def\lognumber#1{\def\thelognumber{#1}}
\def\volumenumber#1{\def\thevolumenumber{#1}}
\def\volumeyear#1{\def\thevolumeyear{#1}}
\def\papernumber#1{\def\thepapernumber{#1}}
\def\pagenumbers#1#2{\def\startpage{#1}\def\finishpage{#2}}
\def\published#1{\def\publishdate{#1}}

\def\received#1{\def\receiveddate{#1}}
\def\revised#1{\def\reviseddate{#1}}
\def\accepted#1{\def\accepteddate{#1}}

\long\def\asciiabstract#1{\long\def\theasciiabstract{#1}}

%  initialise

\let\\\par\let\thelognumber\relax\let\thevolumenumber\relax
\let\thepapernumber\relax\let\thevolumeyear\relax\let\startpage\relax
\let\finishpage\relax\let\publishdate\relax\let\receiveddate\relax
\let\reviseddate\relax\let\accepteddate\relax\let\theasciititle\relax
\let\theasciiauthors\relax
\let\theasciiabstract\relax

\let\theasciiemail\relax

%%%% fonts for AGT logo:

\ifplaintex
\font\logobig=cmssbx10 scaled 3836
\font\logomed=cmssbx10 scaled 2557
\else
\font\logobig=cmssbx10 scaled 4200
\font\logomed=cmssbx10 scaled 2800
\fi

\long\def\makeagttitle{   %%% start of definition of \makeagttitle
\count0=\startpage
\agt\hfill      %   Journal title (top left) 
%   logo (top right)
\hbox to 45truept{\vbox to 0pt{\vglue -13truept{\logomed A\kern -.37em{\logobig 
T}\kern -.38em G}\vss}\hss}
\break
{\small Volume \thevolumenumber\ (\thevolumeyear)
\startpage--\finishpage\nl
Published: \publishdate}

\vglue .25truein

% title
{\parskip=0pt\leftskip 0pt plus
1fil\def\\{\par\smallskip}{\Large\bf\thetitle}\par\medskip} \vglue
0.05truein

% authors :
%
{\parskip=0pt\leftskip 0pt plus 1fil\def\\{\par}{\sc\theauthors}
\par\medskip}%
 
\vglue 0.03truein 

%  abstract and classification numbers:

{\small\leftskip 25truept\rightskip 25truept{\bf Abstract}\stdspace\theabstract

{\bf AMS Classification}\stdspace\theprimaryclass
\ifx\thesecondaryclass\relax\else; \thesecondaryclass\fi\par
{\bf Keywords}\stdspace \thekeywords\par}\vglue 7truept

}   %%%% end of definition of \makeagttitle

\ifplaintex
%  get print centerpage:
\hoffset 14truemm
\voffset 31truemm
%  fonts for headline and footline
\font\phead=cmsl9 scaled 950
\font\pnum=cmbx10 scaled 913
\font\pfoot=cmsl9 scaled 950
%  headline and footline
\headline{\vbox to 0pt{\vskip -4.5mm\line{\small\phead\ifnum
\count0=\startpage ISSN 1472-2739 (on-line) 1472-2747 (printed)
\hfill {\pnum\folio}\else\ifodd\count0\def\\{ }% 
\ifx\theshorttitle\relax\thetitle\else\theshorttitle\fi\hfill{\pnum\folio}
\else\def\\{ and }{\pnum\folio}\hfill\ifx\theshortauthors\relax\theauthors
\else\theshortauthors\fi\fi\fi}\vss}}
\footline{\vbox to 0pt{\vglue 0mm\line{\small\pfoot\ifnum\count0=\startpage
\copyright\ \gtp\hfill\else
\agt, Volume \thevolumenumber\ (\thevolumeyear)\hfill\fi}\vss}}
\else
%  get print centerpage:
\headsep 23pt
\footskip 35pt
\hoffset -4truemm
\voffset 12.5truemm
%  fonts for headline and footline
\font\lhead=cmsl9 scaled 1050
\font\lnum=cmbx10 
\font\lfoot=cmsl9 scaled 1050
\makeatletter
%  headline and footline
\def\@oddhead{{\small\lhead\ifnum\count0=\startpage ISSN 1472-2739 
(on-line) 1472-2747 (printed)\hfill {\lnum\number\count0}\else\ifodd\count0
\def\\{ }\ifx\theshorttitle\relax \thetitle \else\theshorttitle\fi\hfill
{\lnum\number\count0}\else\def\\{ and }{\lnum\number\count0}
\hfill\ifx\theshortauthors\relax 
\theauthors\else\theshortauthors\fi\fi\fi}}\def\@evenhead{\@oddhead}
\def\@oddfoot{\small\lfoot\ifnum\count0=\startpage\copyright\ \gtp\hfill\else
\agt, Volume \thevolumenumber\ (\thevolumeyear)\hfill\fi}
\def\@evenfoot{\@oddfoot}
\makeatother
\fi
%  force \makeagttitle
\let\maketitlepage\makeagttitle
\let\makeshorttitle\maketitlepage
\let\maketitle\maketitlepage

   %%%comment out to create xxx header file

\newwrite\gtoutfile
\long\gdef\makeheadfile{  %%% start of definition of \makeheadfile
{\def\\{, }\def\s{ }
\immediate\openout\gtoutfile head.xxx
\immediate\write\gtoutfile{To: math@arxiv.org}
\immediate\write\gtoutfile{Subject: put OR rep NNNNN:ppppp}
\immediate\write\gtoutfile{--text follows this line--}
\immediate\write\gtoutfile{Proxy-for: \ifx\theasciiauthors\relax
\theauthors\else\theasciiauthors\fi\s<\ifx\theasciiemail\relax\theemail\else\theasciiemail\fi>}
\immediate\write\gtoutfile{\noexpand\\}
\immediate\write\gtoutfile{Authors: \ifx\theasciiauthors\relax
\theauthors\else\theasciiauthors\fi}
{\def\\{ }\immediate\write\gtoutfile{Title: \ifx\theasciititle\relax
\thetitle\else\theasciititle\fi}}
\immediate\write\gtoutfile{Subj-class: GT or SG, GR etc}
\immediate\write\gtoutfile{MSC-class: \theprimaryclass\ifx\thesecondaryclass\relax\else, \thesecondaryclass\fi}
\immediate\write\gtoutfile{Journal-ref: Algebraic and Geometric Topology \thevolumenumber\s
(\thevolumeyear) \startpage-\finishpage}
\immediate\write\gtoutfile{Comments: Published by Algebraic and
Geometric Topology at}
\immediate\write\gtoutfile{\s\s\s  http://www.maths.warwick.ac.uk/agt/AGTVol\thevolumenumber/agt-\thevolumenumber-\thepapernumber.abs.html}
\immediate\write\gtoutfile{\noexpand\\}
\immediate\write\gtoutfile{}
\ifx\theasciiabstract\relax
\immediate\write\gtoutfile{\theabstract}\else
\immediate\write\gtoutfile{\theasciiabstract}\fi
\immediate\write\gtoutfile{}
\immediate\write\gtoutfile{\noexpand\\}
\immediate\write\gtoutfile{}
\immediate\closeout\gtoutfile}}  %%% end of definition of \makeheadfile

\def\maketitlepage{\makeagttitle\makeheadfile}
\let\makeshorttitle\maketitlepage
\let\maketitle\maketitlepage

%%%
%%%  This is agtout.tex.  
%%%
%%%  This the version of  gtoutput.tex  intended to finish formatting
%%%  papers published in Algebriac & Geometric Topology and stored in the
%%%  arXiv.   All versions of  gtoutput.tex  are copyright 
%%%  GT Publications and are to be used _only_ for formatting
%%%  the officially published version of ABT or G&T papers.
%%%
%%%
%%%                                             Colin Rourke  27.102000
%%%
%%%  To create header file  head.xxx  comment out the first \endinput

%  test for latex or plain tex
\def\ifplaintex{\expandafter\ifx\csname documentclass\endcsname\relax}

\def\gtp{{\mathsurround=0pt\it $\cal G\mskip-2mu$eometry \&\ 
$\cal T\!\!$opology $\cal P\!$ublications}}  % GT publications

\def\recd{{\small Received:\qua\receiveddate\ifx\reviseddate\relax
\else\qquad Revised:\qua\reviseddate\fi\par}} 

%  define the various new ingredients of the title page and the data
%  output files

\def\lognumber#1{\def\thelognumber{#1}}
\def\volumenumber#1{\def\thevolumenumber{#1}}
\def\volumeyear#1{\def\thevolumeyear{#1}}
\def\papernumber#1{\def\thepapernumber{#1}}
\def\pagenumbers#1#2{\def\startpage{#1}\def\finishpage{#2}}
\def\published#1{\def\publishdate{#1}}

\def\received#1{\def\receiveddate{#1}}
\def\revised#1{\def\reviseddate{#1}}
\def\accepted#1{\def\accepteddate{#1}}

\long\def\asciiabstract#1{\long\def\theasciiabstract{#1}}

%  initialise

\let\\\par\let\thelognumber\relax\let\thevolumenumber\relax
\let\thepapernumber\relax\let\thevolumeyear\relax\let\startpage\relax
\let\finishpage\relax\let\publishdate\relax\let\receiveddate\relax
\let\reviseddate\relax\let\accepteddate\relax\let\theasciititle\relax
\let\theasciiauthors\relax
\let\theasciiabstract\relax

\let\theasciiemail\relax

%%%% fonts for AGT logo:

\ifplaintex
\font\logobig=cmssbx10 scaled 3836
\font\logomed=cmssbx10 scaled 2557
\else
\font\logobig=cmssbx10 scaled 4200
\font\logomed=cmssbx10 scaled 2800
\fi

\long\def\makeagttitle{   %%% start of definition of \makeagttitle
\count0=\startpage
\agt\hfill      %   Journal title (top left) 
%   logo (top right)
\hbox to 45truept{\vbox to 0pt{\vglue -13truept{\logomed A\kern -.37em{\logobig 
T}\kern -.38em G}\vss}\hss}
\break
{\small Volume \thevolumenumber\ (\thevolumeyear)
\startpage--\finishpage\nl
Published: \publishdate}

\vglue .25truein

% title
{\parskip=0pt\leftskip 0pt plus
1fil\def\\{\par\smallskip}{\Large\bf\thetitle}\par\medskip} \vglue
0.05truein

% authors :
%
{\parskip=0pt\leftskip 0pt plus 1fil\def\\{\par}{\sc\theauthors}
\par\medskip}%
 
\vglue 0.03truein 

%  abstract and classification numbers:

{\small\leftskip 25truept\rightskip 25truept{\bf Abstract}\stdspace\theabstract

{\bf AMS Classification}\stdspace\theprimaryclass
\ifx\thesecondaryclass\relax\else; \thesecondaryclass\fi\par
{\bf Keywords}\stdspace \thekeywords\par}\vglue 7truept

}   %%%% end of definition of \makeagttitle

\ifplaintex
%  get print centerpage:
\hoffset 14truemm
\voffset 31truemm
%  fonts for headline and footline
\font\phead=cmsl9 scaled 950
\font\pnum=cmbx10 scaled 913
\font\pfoot=cmsl9 scaled 950
%  headline and footline
\headline{\vbox to 0pt{\vskip -4.5mm\line{\small\phead\ifnum
\count0=\startpage ISSN 1472-2739 (on-line) 1472-2747 (printed)
\hfill {\pnum\folio}\else\ifodd\count0\def\\{ }% 
\ifx\theshorttitle\relax\thetitle\else\theshorttitle\fi\hfill{\pnum\folio}
\else\def\\{ and }{\pnum\folio}\hfill\ifx\theshortauthors\relax\theauthors
\else\theshortauthors\fi\fi\fi}\vss}}
\footline{\vbox to 0pt{\vglue 0mm\line{\small\pfoot\ifnum\count0=\startpage
\copyright\ \gtp\hfill\else
\agt, Volume \thevolumenumber\ (\thevolumeyear)\hfill\fi}\vss}}
\else
%  get print centerpage:
\headsep 23pt
\footskip 35pt
\hoffset -4truemm
\voffset 12.5truemm
%  fonts for headline and footline
\font\lhead=cmsl9 scaled 1050
\font\lnum=cmbx10 
\font\lfoot=cmsl9 scaled 1050
\makeatletter
%  headline and footline
\def\@oddhead{{\small\lhead\ifnum\count0=\startpage ISSN 1472-2739 
(on-line) 1472-2747 (printed)\hfill {\lnum\number\count0}\else\ifodd\count0
\def\\{ }\ifx\theshorttitle\relax \thetitle \else\theshorttitle\fi\hfill
{\lnum\number\count0}\else\def\\{ and }{\lnum\number\count0}
\hfill\ifx\theshortauthors\relax 
\theauthors\else\theshortauthors\fi\fi\fi}}\def\@evenhead{\@oddhead}
\def\@oddfoot{\small\lfoot\ifnum\count0=\startpage\copyright\ \gtp\hfill\else
\agt, Volume \thevolumenumber\ (\thevolumeyear)\hfill\fi}
\def\@evenfoot{\@oddfoot}
\makeatother
\fi
%  force \makeagttitle
\let\maketitlepage\makeagttitle
\let\makeshorttitle\maketitlepage
\let\maketitle\maketitlepage

   %%%comment out to create xxx header file

\newwrite\gtoutfile
\long\gdef\makeheadfile{  %%% start of definition of \makeheadfile
{\def\\{, }\def\s{ }
\immediate\openout\gtoutfile head.xxx
\immediate\write\gtoutfile{To: math@arxiv.org}
\immediate\write\gtoutfile{Subject: put OR rep NNNNN:ppppp}
\immediate\write\gtoutfile{--text follows this line--}
\immediate\write\gtoutfile{Proxy-for: \ifx\theasciiauthors\relax
\theauthors\else\theasciiauthors\fi\s<\ifx\theasciiemail\relax\theemail\else\theasciiemail\fi>}
\immediate\write\gtoutfile{\noexpand\\}
\immediate\write\gtoutfile{Authors: \ifx\theasciiauthors\relax
\theauthors\else\theasciiauthors\fi}
{\def\\{ }\immediate\write\gtoutfile{Title: \ifx\theasciititle\relax
\thetitle\else\theasciititle\fi}}
\immediate\write\gtoutfile{Subj-class: GT or SG, GR etc}
\immediate\write\gtoutfile{MSC-class: \theprimaryclass\ifx\thesecondaryclass\relax\else, \thesecondaryclass\fi}
\immediate\write\gtoutfile{Journal-ref: Algebraic and Geometric Topology \thevolumenumber\s
(\thevolumeyear) \startpage-\finishpage}
\immediate\write\gtoutfile{Comments: Published by Algebraic and
Geometric Topology at}
\immediate\write\gtoutfile{\s\s\s  http://www.maths.warwick.ac.uk/agt/AGTVol\thevolumenumber/agt-\thevolumenumber-\thepapernumber.abs.html}
\immediate\write\gtoutfile{\noexpand\\}
\immediate\write\gtoutfile{}
\ifx\theasciiabstract\relax
\immediate\write\gtoutfile{\theabstract}\else
\immediate\write\gtoutfile{\theasciiabstract}\fi
\immediate\write\gtoutfile{}
\immediate\write\gtoutfile{\noexpand\\}
\immediate\write\gtoutfile{}
\immediate\closeout\gtoutfile}}  %%% end of definition of \makeheadfile

\def\maketitlepage{\makeagttitle\makeheadfile}
\let\makeshorttitle\maketitlepage
\let\maketitle\maketitlepage

\lognumber{4}
\volumenumber{1}
\volumeyear{2001}
\papernumber{4}
\published{27 January 2001}
\pagenumbers{57}{71}
\received{11 December 2000}
\revised{12 January 2001}
\accepted{12 January 2001}

%%%% user macros
\def\ni\noindent
\def\sbs{\subset}

%\def\as{\operatorname{asdim}}

%\def l{\operatorname{vol}}
%\def\dist{\operatorname{dist}}

%
%   End of user macros
%
%
%%%%%%%%%%%%          Main text starts here        %%%%%%%%%
%
%                       References
%
\reflist

\key{D}
{\bf M. Davis}, 
{\it Groups generated by reflections and aspherical manifolds not covered
by Euclidean space},
 Ann. Math.,
 (2) 117,
(1983),
 293-325.

\key{D2}
{\bf M. Davis}, 
{\it Coxeter groups and aspherical manifolds},
Lecture Notes in Math.
 1051
(1984),
 197-221.

\key{D-H}
{\bf M. Davis and J.-C. Hausmann}, 
{\it Aspherical manifolds without smooth or PL structure,}
 Lecture Notes in Math.,
 (2) 1370 (1989), 135-142.

\key{Dr1}
{\bf A. Dranishnikov,}
 {\it Asymptotic topology}
  Russian Mathematical Surveys,
  55 (6),
  (2000).

\key{Dr2}
 {\bf A. Dranishnikov,}
 {\it On large scale properties of manifolds,}
  Preprint, 1999.\nl
  {\tt arxiv:math.GT/9912062}

\key{Dr3}
 {\bf A. Dranishnikov},
 {\it Hypereuclidean manifolds and the Novikov Conjecture,}\nl
 {Preprint of MPI, Bonn,}
  2000.

\key{DJ}
 {\bf A. Dranishnikov and T. Januszkiewicz},
 {\it Every Coxeter group acts amenably on a compact space,}
  Preprint, to appear in Topology Proceedings.\nl
  {\tt arxiv:math.GT/9911245}

\key{G1}   {\bf M. Gromov,}
 {\it Asymptotic invariants of infinite groups}, 
  Cambridge University Press,
 Geometric Group Theory, vol 2
  1993. 

\key{G2}
 {\bf M. Gromov},
 {\it Spaces and questions},
  Preprint,
  1999.

\key{M}
 {\bf G. Mess},
 {\it Examples of Poincare duality groups},
  Proc. Amer. Math. Soc.
  110, no 4,
  (1990), 1145-1146.

\key{R}
 {\bf J. Roe},
{\it Coarse cohomology and index theory for complete Riemannian 
manifolds}, 
  Memoirs Amer. Math. Soc. No. 497,
  1993.

\key{S}
 {\bf J.-P. Serre},
{\it Trees},
  Springer-Verlag,
  1980.

\key{Y}
 {\bf G. Yu},
 {\it The Novikov conjecture for groups with finite asymptotic 
dimension},
  Ann. of Math.,
  147, no 2,
  (1998), 325-355.

\endreflist
%
%                      Header material (title etc)  

\title{On asymptotic dimension of groups}                    
\authors{G. Bell and A. Dranishnikov}                  
\address{University of Florida, Department of Mathematics,\\
PO Box 118105,
358 Little Hall,\\Gainesville, FL 32611-8105, USA}                  
\email{dranish@math.ufl.edu, bell@math.ufl.edu}

\abstract 
We prove a version of the countable union theorem for asymptotic dimension
and we apply it to groups acting on asymptotically finite dimensional metric spaces.
As a consequence we obtain the following finite dimensionality theorems.

A)\qua An amalgamated product of asymptotically finite dimensional groups has finite
asymptotic dimension: $asdim A\ast_CB<\infty$.

B)\qua  Suppose that $G'$ is an HNN extension of a group $G$ with 
$asdim G<\infty$. 
Then $asdim G'<\infty$.

C)\qua  Suppose that $\Gamma$ is Davis' group constructed from a group $\pi$ with
$asdim\pi<\infty$. Then 
$asdim\Gamma<\infty$.

\endabstract

\asciiabstract{
We prove a version of the countable union theorem for asymptotic
dimension and we apply it to groups acting on asymptotically finite
dimensional metric spaces.  As a consequence we obtain the following
finite dimensionality theorems.

A) An amalgamated product of asymptotically finite dimensional groups has 
finite asymptotic dimension: asdim A *_C B < infinity.

B) Suppose that G' is an HNN extension of a group G with asdim G < infinity. 
Then asdim G'< infinity.

C) Suppose that \Gamma is Davis' group constructed from a group \pi with
asdim\pi < infinity. Then asdim\Gamma < infinity.}

\primaryclass{20H15}                
\secondaryclass{20E34, 20F69}              
\keywords{Asymptotic dimension, amalgamated product, HNN extension}                    

\maketitle

\def\\{\par}
\section{Introduction}

The notion of the asymptotic dimension was introduced by Gromov [\ref{G1}]
as an asymptotic analog of Ostrand's characterization of covering
dimension. Two sets $U_1$, $U_2$ in a metric space are called $d$-disjoint
if they are at least $d$-apart, i.e. $\inf\{dist(x_1,x_2)\mid x_1\in U_1,
x_2\in U_2\}\ge d$.
A metric space $X$ has  {\it asymptotic dimension} 
$asdim X\le n$ if for
an arbitrarily large number $d$ one can find $n+1$ uniformly bounded families
${\cal U}^0,\dots,{\cal U}^n$ of $d$-disjoint sets in $X$ such that the 
union $\cup_i{\cal U}^i$ is
a cover of $X$. A generating set $S$ in a group $\Gamma$ defines the 
{\it word metric}
on $\Gamma$ by the following rule: 
$d_S(x,y)$ is the minimal length of a presentation of the element
$x^{-1}y\in\Gamma$ in the alphabet $S$.
Gromov applied the notion of asymptotic dimension
to studying asymptotic invariants of discrete groups. 
It follows from the definition that the asymptotic dimension
$asdim(\Gamma, d_S)$ of a finitely generated group does not depend on 
the choice of the finite
generating set $S$. Thus, $asdim\Gamma$ is an asymptotic invariant for finitely 
generated groups. Gromov proved [\ref{G1}] that $asdim\Gamma<\infty$ for hyperbolic groups 
$\Gamma$. The corresponding question about nonpositively curved (or CAT(0)) groups 
remains
open. In the case of Coxeter groups it was answered in [\ref{DJ}].

In [\ref{Y}] G. Yu proved a series of conjectures, including the famous
Novikov Higher Signature conjecture, for groups $\Gamma$ with
$asdim\Gamma<\infty$. Thus, the problem of determining the 
asymptotic finite dimensionality
of certain discrete groups became very important. 
In fact, until the recent example of
Gromov [\ref{G2}] it was unknown whether all
finitely presented groups satisfy the inequality $asdim\Gamma<\infty$.
In view of this, it is natural to ask whether the property of 
asymptotic finite dimensionality
is preserved under the standard constructions with groups.
Clearly, the answer is positive for the direct product of two groups.
It is less clear, but still is not difficult to see that a semidirect product
of asymptotically finite dimensional groups has a finite asymptotic dimension.
The same question about the free product does not seem clear at all. 
In this paper
we show that the asymptotic finite dimensionality is preserved by 
the free product,
by the amalgamated free product and by the HNN extension.

One of the motivations  for this paper was to prove that Davis' construction 
preserves asymptotic finite dimensionality. Given a group $\pi$ with a finite 
classifying space $B\pi$,  Davis found a canonical construction, based on
 Coxeter groups, of a group $\Gamma$ with $B\Gamma$ a closed manifold such that
$\pi$ is a retract of $\Gamma$ (see [\ref{D}],[\ref{D2}],[\ref{D-H}],[\ref{M}]). We prove here that if $asdim\pi<\infty$, then 
$asdim\Gamma<\infty$. This theorem together with the result
of the second author [\ref{Dr3}] (see also [\ref{Dr2}]) about the hypereuclideanness of 
asymptotically finite dimensional manifolds allows one to get a shorter and more 
elementary proof of the Novikov Conjecture for groups $\Gamma$ with 
$asdim\Gamma<\infty$.

We note that the asymptotic dimension $asdim$ is a coarse invariant, i.e. it is an 
invariant of the {\it coarse category} introduced in [\ref{R}]. We recall that the objects 
in the coarse category are metric spaces and morphisms are coarsely proper 
and coarsely uniform (not necessarily continuous) maps. A map $f:X\to Y$ 
between metric spaces is called {\it coarsely proper} if the preimage 
$f^{-1}(B_r(y))$ of
every ball in $Y$ is a bounded set in $X$. A map $f:X\to Y$ is called 
{\it coarsely uniform} if
there is a function $\rho:{\bf R}_+\to{\bf R}_+$, tending to infinity, such that
$d_Y(f(x),f(y))\le\rho(d(x,y))$ for all $x,y\in Y$.  We note 
that every object in the
coarse category is isomorphic to a discrete metric space.

There is an analogy between the standard (local) topology and the
asymptotic topology which is outlined in [\ref{Dr1}]. That analogy is not 
always direct.
Thus, in Section 2 we prove the following finite union theorem 
for asymptotic dimension $asdim X\cup Y\le\max\{asdim X,asdim Y\}$ whereas 
the classical Menger-Urysohn theorem states: $\dim X\cup Y\le\dim X+\dim Y+1$.
Also the Countable Union Theorem in the classical dimension theory cannot have
a straightforward analog, since all interesting objects in
the coarse category are countable unions of points
but not all of them are asymptotically 0-dimensional. In Section 2 we formulated
a countable union theorem for asymptotic dimension
which we found useful for applications to the case of discrete groups.

The second author was partially supported by NSF grant DMS-9971709.

\section{Countable union theorem}

\rk{Definition} A family of metric spaces $\{F_{\alpha}\}$ satisfies
the inequality $asdim F_{\alpha}\le n$ {\it uniformly} if for arbitrarily 
large 
$d>0$ there are $R$ and $R$-bounded $d$-disjoint families 
${\cal U}^0_{\alpha}\dots {\cal U}_{\alpha}^n$ of subsets of $F_{\alpha}$ 
such that the
union $\cup_i{\cal U}^i_{\alpha}$ is a cover of $F_{\alpha}$.

A typical example of such family is when all $F_{\alpha}$ are isometric to
a space $F$ with $asdim F\le n$. 

A discrete metric space $X$ has {\it bounded geometry} if for every $R$ there
is a constant $c=c(R)$ such that every $R$-ball $B_R(x)$ in $X$ 
contains at most $c$ points.
\proclaim{Proposition 1}
Let $f_{\alpha}:F_{\alpha}\to X$ be a family of 1-Lipschitz injective maps 
to a discrete metric space
of bounded geometry with $asdim X\le n$. Then $asdim F_{\alpha}\le n$ uniformly.
\endproclaim
\prf
For a metric space $A$ we define its $d$-components as the classes under
the following equivalence relation. Two points $a,a'\in A$ are equivalent if
there is a chain of points $a_0, a_1,\dots,a_k$ with $a_0=a$, $a_k=a'$ and with
$d(a_i,a_{i+1})\le d$ for all $i<k$. We note that the $d$-components are 
more than $d$ apart and also note that the diameter of each $d$-component is
less than or equal to $d|A|$, where $|A|$ is the number of points in $A$.

Let $d$ be given. Then there are $R$-bounded $d$-disjoint families
${\cal V}^0,\dots,{\cal V}^n$ covering $X$. For every $V\in{\cal V}^i$ and 
every $\alpha$
we present the set $f_{\alpha}^{-1}(V)$ as the union of $d$-components:
$f_{\alpha}^{-1}(V)=\cup C_{\alpha}^j(V)$. Note that the diameter of every 
$d$-component
is $\le dc(R)$ where the function $c$ is taken from the bounded 
geometry condition on $X$.
 We take ${\cal U}^i_{\alpha}=\{C^j_{\alpha}(V)\mid V\in {\cal V}^i\}$. 
\endprf

\proclaim{Theorem 1}
Assume that $X=\cup_{\alpha}F_{\alpha}$ and $asdim F_{\alpha}\le n$ uniformly.
Suppose that for any $r$ there exists $Y_r\subset X$ with $asdim Y_r\le n$ and
such that the family $\{F_{\alpha}\setminus Y_r\}$ is $r$-disjoint.
Then $asdim X\le n$.
\endproclaim
\proclaim{Finite Union Theorem}
Suppose that a metric space is presented as a union $A\cup B$ of subspaces. 
Then $asdim A\cup B\le\max\{asdim A,asdim B\}$.
\endproclaim
\prf
We apply Theorem 1 to the case when the
 family of subsets consists of $A$ and $B$ and we take $Y_r=B$.\endprf

The proof of Theorem 1 is based on the idea of saturation of one family by the other.
Let $\cal V$ and $\cal{U}$ be two families of subsets of a metric space $X$.

\rk{Definition} For $V\in{\cal V}$ and $d>0$ we denote by $N_d(V;{\cal U})$ the union of $V$
and all elements $U\in{\cal U}$ with 
$d(V,U)=\min\{d(x,y)\mid x\in V,y\in U\}\le d$.
By $d$-saturated union of ${\cal V}$ and ${\cal U}$ we mean the following
family ${\cal V}\cup_d{\cal U}=\{N_d(V;{\cal U})\mid V\in{\cal V}\}
\cup\{U\in{\cal U}\mid d(U,V)>d\ \ for\ \ all\ \ V\in{\cal V}\}$.

Note that this is not a commutative operation. Also note that
$\{\emptyset\}\cup_d{\cal U}={\cal U}$ and 
${\cal V}\cup_d\{\emptyset\}={\cal V}$ for all $d$.

\proclaim{Proposition 2}
Assume that $\cal U$ is $d$-disjoint and $R$-bounded, $R\ge d$. Assume that 
$\cal V$ is
$5R$-disjoint and $D$-bounded. Then ${\cal V}\cup_d{\cal U}$ is $d$-disjoint and
$D+2(d+R)$-bounded.
\endproclaim
\prf
First we note that elements of type $U$ are $d$-disjoint in the saturated union.
The same is true for elements of type $U$ and $N_d(V;{\cal U})$. Now consider elements
$N_d(V;{\cal U})$ and $N_d(V';{\cal U})$. Note that they are contained in the 
$d+R$-neighborhoods of $V$ and $V'$ respectively. Since $V$ and $V'$ are
$5R$-disjoint, and $R\ge d$, the neighborhoods will be $d$-disjoint.

Clearly, $diam N_d(V;{\cal U})\le diam V+2(d+R)\le D+2(d+R)$.
\endprf

{\bf Proof of Theorem 1}\qua
Let $d$ be given. Consider $R$ and families 
${\cal U}^0_{\alpha}\dots{\cal U}^n_{\alpha}$
from the definition of the uniform inequality $asdim F_{\alpha}\le n$. 
We may assume that $R>d$. We take $r=5R$ and consider $Y_r$ satisfying the 
conditions of the Theorem.
Consider $r$-disjoint $D$-bounded families ${\cal V}^0,\dots,{\cal V}^k$ 
from the definition of $asdim Y_r\le k$. Let $\bar{\cal U}^i_{\alpha}$ 
be the restriction of
${\cal U}^i_{\alpha}$ over $F_{\alpha}\setminus Y_r$, i.e.
$\bar{\cal U}^i_{\alpha}=\{U\setminus Y_r\mid U\in{\cal U}^i_{\alpha}\}$.
Let $\bar{\cal U}^i=\cup_{\alpha}\bar{\cal U}^i_{\alpha}$. Note that the family
$\bar{\cal U}^i$ is $d$-disjoint and $R$-bounded.
For every $i$ we define ${\cal W}^i={\cal V}^i\cup_d\bar{\cal U}^i$ . 
By Proposition 2 the family
${\cal W}^i$ is $d$-disjoint and uniformly bounded.
Clearly $\cup_i{\cal W}^i$ covers $X$.\endprf

\section{Groups acting on finite dimensional spaces}

A norm on a group $A$ is a map $\|\ \|:A\to{\bf Z}_+$ such that
$\|ab\|\le\|a\|+\|b\|$ and $\|x\|=0$ if and only if $x$ is the unit in $A$.
A set of generators $S\subset A$ defines the norm $\|x\|_S$ as the minimal
length of a presentation of $x$ in terms of $S$. A norm on a group defines 
 a left invariant metric $d$ by $d(x,y)=\|x^{-1}y\|$. 
If $G$ is a finitely generated group and $S$ and $S'$ are two finite generating sets,
then the corresponding metrics $d_S$ and $d_{S'}$ define coarsely equivalent
metric spaces $(G,d_S)$ and $(G,d_{S'})$. In particular, 
$asdim(G,d_S)=asdim(G,d_{S'})$, and we can speak about the asymptotic dimension
$asdim G$ of a finitely generated group $G$.

 Assume that a group $\Gamma$ acts on a metric space
$X$. For every $R>0$ we define the $R$-stabilizer $W_R(x_0)$ of a point 
$x_0\in X$ as the set of all $g\in\Gamma$ with $g(x_0)\in B_R(x_0)$.
Here $B_R(x)$  denotes the closed ball of radius $R$ centered at $x$.

\proclaim{Theorem 2}
Assume that a finitely generated group $\Gamma$ acts by isometries on
a metric space $X$ with a base point $x_0$ and with $asdim X\le k$.
Suppose that $asdim W_R(x_0)\le n$ for all $R$.
Then $asdim\Gamma\le (n+1)(k+1)-1$.
\endproclaim
\prf
We define a map $\pi:\Gamma\to X$ by the formula $\pi(g)=g(x_0)$.
Then $W_R(x_0)=\pi^{-1}(B_r(x_0))$.
Let $\lambda=\max\{d_X(s(x_0),x_0)\mid s\in S\}$. We show now that $\pi$ is
$\lambda$-Lipschitz. Since the metric $d_S$ on $\gamma$ is induced
from the geodesic metric on the Cayley graph, it suffices to check that
$d_X(\pi(g),\pi(g'))\le\lambda$ for all $g,g'\in\Gamma$ with $d_S(g,g')=1$.
Without loss of generality we assume that $g'=gs$ where $s\in S$.
Then $d_X(\pi(g),\pi(g'))=d_X(g(x_0),gs(x_0))=d_X(x_0,s(x_0))\le\lambda$.

Note that $\gamma B_R(x)=B_R(\gamma(x))$ and $\gamma(\pi^{-1}(B_R(x)))=
\pi^{-1}(B_R(\gamma(x)))$ for all $\gamma\in\Gamma$, $x\in X$ and all $R$.

Given $r>0$, there are $\lambda r$-disjoint,  $R$-bounded families
${\cal F}^0,\dots,{\cal F}^k$ on the orbit $\Gamma x_0$. Let 
 ${\cal V}^0,\dots,{\cal V}^n$ on $W_{2R}(x_0)$
 be $r$-disjoint uniformly bounded families given by the definition of
the inequality $asdim W_R(x_0)\le n$. For every element $F\in{\cal F}^i$ we 
choose an element $g_F\in\Gamma$ such that $g_F(x_0)\in F$. We define
$(k+1)(n+1)$ families of subsets of $\Gamma$ as follows:
$$
{\cal W}^{ij}=\{g_F(C)\cap\pi^{-1}(F)\mid F\in{\cal F}^i, C\in{\cal V}^j\}
$$ 
Since  multiplication by $g_F$ from the left is an isometry, every 
two distinct sets $g_F(C)$ and $g_F(C')$
are $r$-disjoint. Note that $\pi(g_F(C)\cap\pi^{-1}(F))$ and 
$\pi(g_{F'}(C)\cap\pi^{-1}(F'))$ are $\lambda r$-disjoint
for $F\ne F'$. Since $\pi$ is $\lambda$-Lipschitz, the sets 
$g_F(C)\cap\pi^{-1}(F)$ and
$g_{F'}(C')\cap\pi^{-1}(F')$ are $r$-disjoint. 
The families ${\cal W}^{ij}$ are uniformly bounded,
since the families ${\cal V}^j$ are, and  multiplication by $g$ from the left is an isometry on $\Gamma$.
We check that the union of the families ${\cal W}^{ij}$ forms a cover of
$\Gamma$. Let $g\in\Gamma$ and let $\pi(g)=F$, i.e. $g(x_0)\in F$. 
Since $diam F\le R$, $x_0\in g_F^{-1}(F)\le R$ and $g_F^{-1}$ acts as 
an isometry, we have $g^{-1}_F(F)\subset B_R(x_0)$. Therefore, 
$g_F^{-1}g(x_0)\in B_R(x_0)$, i. e. $g^{-1}_Fg\in W_R(x_0)$.
Hence $g_F^{-1}g$ lies in some set $C\in{\cal V}^j$ for some $j$.
Therefore $g\in g_F(C)$. Thus, $g\in g_F(C)\cap\pi^{-1}(F)$.
\endprf

\proclaim{Theorem 3}
Let $\phi:G\to H$ be an epimorphism of a finitely generated group $G$ 
with kernel $ker\phi=K$.
Assume that $asdim K\le k$ and $asdim H\le n$. Then $asdim G\le (n+1)(k+1)-1$.
\endproclaim
\prf 
The group $G$ acts on $H$ by the rule $g(h)=\phi(g)h$.
This is an action by isometries for every left invariant metric on $H$.
Let $S$ be a finite generating set for $G$. We consider the metric on $H$ 
defined by the set $\phi(S)$. Below we prove that the 
$R$-stabilizer of the identity $W_R(e)$ coincides with  $N_R(K)$, the 
$R$-neighborhood of $K$ in $G$. Since $N_R(K)$ is coarsely isomorphic to
$K$, we have the inequality $asdim W_R(e)\le k$.

Let $g\in W_R(e)$, then $\|\phi(g)\|\le R$. Therefore, there is a sequence
$i_1,\dots,i_k$ with $k\le R$ such that $\phi(g)=\bar s_{i_1}\dots\bar s_{i_k}$
where $\bar s=\phi(s)$, $s\in S$. Let $u=s_{i_1}\dots s_{i_k}$.
Then $d_S(g,gu^{-1})\le R$ and hence, $d_S(g,K)\le R$. In the opposite direction,
if $d_S(g,K)\le R$, then $d(g,z)\le R$ for some $z\in K$. Hence 
$d_{\phi(S)}(\phi(g),e)\le R$.

We apply Theorem 2 to complete the proof.\endprf

\rk{Remark} The estimate $(n+1)(k+1)-1$ in Theorems 2 and 3 is far from being
sharp. Since in this paper we are interested in finite dimensionality only,
we are not trying to give an exact estimate which is  $n+k$. Besides, it
would be difficult to get an exact estimate just working with covers.
Even for proving the inequality
$$
asdim\Gamma_1\times\Gamma_2\le asdim\Gamma_1+asdim\Gamma_2
$$
it is better to use a different approach to $asdim$ (see [\ref{DJ}]).

\section{Free and amalgamated products}

Let $\{A_i,\|\ \|_i\}$ be a sequence of groups with norms. Then these norms
generate a norm on the free product $\ast A_i$ as follows. Let
$x_{i_1}x_{i_2}\dots x_{i_m}$ be the reduced presentation of $x\in\ast A_i$,
where $x_{i_k}\in A_{i_k}$. We  denote by $l(x)=m$ the length of
the reduced presentation of $x$ and we define 
$\|x\|=\|x_{i_1}\|_{i_1}+\dots+\|x_{i_m}\|_{i_m}$.

\proclaim{Theorem 4}
Let $\{A_i,\|\ \|_i\}$ be a sequence of groups 
satisfying $asdim A_i\le n$ uniformly
 and let $\|\ \|$ be the norm on the free product
$\ast A_i$ generated by the norms $\|\ \|_i$. 
Then $asdim(\ast A_i,\|\ \|)\le 2n+1$.
\endproclaim
\prf
First we note that the uniform property $asdim A_i\le n$ and Theorem 1
applied with $Y_r=B_r(e)$, 
the $r$-ball in $\ast A_i$ centered at the unit $e$,
imply that $asdim\cup A_i\le n$.

We let $G$ denote $\ast A_i$. Then we consider a tree $T$ with vertices left
cosets $xA_j$ in $G$. Two vertices $xA_i$ and $yA_j$ are joined by an edge 
if and only if there is an element $z\in G$ such that $xA_i=zA_i$ and
$yA_j=zA_j$ and $i\ne j$. The multiplication by elements of $G$ from the left
defines an action of $G$ on $T$. We note that the $m$-stabilizer $W_m(A_1)$ 
of the vertex $A_1$ is the union of all possible products $A_{i_1}\dots A_{i_l}A_1$ 
of the length
$\le m+1$, where $i_k\ne i_{k+1}$ and $i_l\ne 1$. 
Let $P_m=\{x\in \ast A_i\mid l(x)=m\}$ and let
$P^k_m=\{x\in P_m\mid x=x_{i_1}\dots x_{i_m}, x_{i_m}\notin A_k\}$.
Put $R_m=W_m(A_1)\setminus W_{m-1}(A_1)$.
Then $R_m\subset P_{m+1}$

By induction on $m$ we show that
$asdim P_m\le n$. This statement holds true when $m=0$, since $P_0=\{e\}$.
Assume that it holds for $P_{m-1}$.
We note that $P_m=\cup_{x\in P^i_{m-1}}xA_i$.
Since multiplication from the left is an isometry, the hypothesis
of the theorem implies that the inequality $asdim xA_i\le n$ holds uniformly.
Given $r$ we consider the set $Y_r=P_{m-1}B_r(e)$
where  $B_r(e)$ is the $r$-ball in $\ast A_i$.
Since $Y_r$ contains $P_{m-1}$
and is contained in $r$-neighborhood of $P_{m-1}$, it is isomorphic in
the coarse category to $P_{m-1}$. Hence by the induction assumption
we have $asdim Y_r\le n$.
We show that the family $xA_i\setminus  Y_r$, $x\in P_{m-1}^i$ is
$r$-disjoint. Assume that $xA_i\ne x'A_j$. This means that $x\ne x'$ if 
$i=j$. If $i\ne j$ the inequality $\|a_i^{-1}x^{-1}x'a_j\|\ge\|a_i^{-1}a_j\|=
\|a_i\|+\|a_j\|$
holds for any choice of $a_i\in A_i$ and $a_j\in A_j$. If $i=j$, the same 
inequality holds,
since $x\ne x'$ and they are of the same length.
If $xa_i\in xA_i\setminus  Y_r$ and $x'a_j\in xA_j\setminus  Y_r$,
then $\|a_i\|,\|a_j\|\ge r$ and hence $ dist(xa_i,x'a_j)\ge 2r$.
Theorem 1 implies that $asdim P_m\le n$. The Finite Union Theorem implies that
$asdim W_m(A_1)\le n$ for all $n$.

It is known that every tree $T$ has $asdim T=1$ (see [\ref{DJ}]).
Thus by Theorem 2 $asdim(\ast A_i,\|\ \|)\le 2n+1$.
\endprf

\proclaim{Corollary}
Let $A_i$, $i=1,\dots, k$, be finitely generated groups with $asdim A_i\le n$.
Then $asdim\ast_{i=1}^kA_i\le 2n+1$.
\endproclaim

\proclaim{Theorem 5} Let $A$ and $B$ be finitely generated groups
with $asdim A\le n$ and $asdim B\le n$
and let $C$ be their common subgroup. Then $asdim A\ast_C B\le 2n+1$.
\endproclaim

We recall that every element $x\in A\ast_C B$ admits a  unique 
normal presentation
$c\bar x_1\dots\bar x_k$ where $c\in C$, $\bar x_i=Cx_i$ are nontrivial 
alternating right cosets of $C$ in $A$ or $B$.
Thus, $x=cx_1\dots x_k$. Let $ dist(\ ,\ )$ be a metric on the group $G=A\ast_CB$.
We assume that this metric is generated by the union of the finite sets of
generators $S=S_A\cup S_B$ of the groups $A$ and $B$.
On the space of the right cosets $C\setminus G$ of a subgroup 
$C$ in $G$
one can define the  metric $\bar d(Cx,Cy)= dist(Cx,Cy)= dist(x,Cy)$. 
The following chain of inequalities implies the triangle inequality for
$\bar d$:
$
 dist(Ca,Cb)\le  dist(a,c'b)=\|a^{-1}c'b\|\le\|a^{-1}cz\|+\|(cz)^{-1}c'b\|
 =dist(a,cz)+ dist(cz,c'b).$
We chose $c$ such that $ dist(a,cz)= dist(a, Cz)=\bar d(Ca,Cz)$ and
$c'$ such that $ dist(cz,c'b)= dist(cz,Cb)=\bar d(Cz,Cb)$.

For every pair of pointed metric spaces $X$ and $Y$ we define a free product
$X\hat\ast Y$ as a metric space whose elements are alternating words formed by
the alphabets $X\setminus \{x_0\}$ and $Y\setminus \{y_0\}$ plus the 
trivial word $x_0=y_0=\tilde e$. 
We define the norm of the trivial word to be zero and for a word 
of type $x_1y_1\dots x_ry_r$ we set 
$\|x_1y_1\dots x_ry_r\|=\Sigma_id_X(x_i,x_0)+d_Y(y_i,y_0)$. If the word starts 
or ends by a different 
type of letter, we consider the corresponding sum.
To define the distance $d(w,w')$ between two words $w$ and $w'$
we cut off their common part $u$ if it is not empty:
$w=uxv$, $w'=ux'v'$ and set $d(w,w')=d(x,x')+\|v\|+\|v'\|$.
If the common part is empty, we define $d(w,w')=\|w\|+\|w'\|$.
Thus, $d(w,\tilde e)=\|w\|$.

\proclaim{Proposition 3}
Let $c\bar x_1\dots\bar x_r$ be the normal presentation of $x\in A\ast_CB$.
Then $\|x\|\ge\Sigma_i\bar d(\bar x_i,C)$.
\endproclaim
\prf
We define a map $\phi:A\ast_CB\to (C\backslash A)\hat\ast(C\backslash B)$ as follows. If 
$c\bar x_1\dots\bar x_r$ is the normal
presentation of $x$, then we set $\phi(x)=\bar x_1\dots\bar x_r$ and
define $\phi(e)=\tilde e$.  
 We verify that $\phi$ is 
1-Lipschitz. Since $A\ast_CB$ is a discrete geodesic metric space space, it suffices to show
that $d(\phi(x),\phi(x\gamma))\le 1$ where $\gamma$ is a generator in $A$
or in $B$. Let $x=cx_1\dots x_r$ be a presentation corresponding to the normal
presentation $c\bar x_1\dots\bar x_r$. Then the normal presentation of $x\gamma$
will be either $c\bar x_1\dots(\overline{x_r\gamma})$ or
$c\bar x_1\dots\bar x_r\bar\gamma$. In the first case,
$d(\phi(x),\phi(x\gamma))=\bar d(\bar x_r,\overline{x_r\gamma})=
 dist(Cx_r,Cx_r\gamma)\le  dist(x_r,x_r\gamma)=1$. In the second case we
have $d(\phi(x),\phi(x\gamma))=\bar d(C,C\gamma)= 
 dist(C,C\gamma)\le  dist(e,\gamma)=1$.

Then $\|x\|= dist(x,e)\ge d(\phi(x),\tilde e)=d(\bar x_1\dots\bar x_r,\tilde e)=
\|x_1\dots x_r\|=\Sigma_i\bar d(\bar x_i,\bar e)$. \endprf

\proclaim{Proposition 4}
Suppose that the subset $(BA)^m=BA\dots BA\subset A\ast_CB$ is supplied with the
induced metric and let $asdim A,asdim B\le n$. Then $asdim(BA)^m\le n$ for all $m$.
\endproclaim
\prf
Let $l(x)$ denote the length of the normal presentation 
$c\bar x_1\dots\bar x_{l(x)}$ of $x$. Define $P_k=\{x\mid l(x)=k\}$,
  $P_k^A=\{x\in P_k\mid x_{l(x)}\in C\backslash A\}$ and 
$P_k^B=\{x\in P_k\mid x_{l(x)}\in C\backslash B\}$. Note that $P_k=P^A_k\cup P^B_k$.
Also we note that $(BA)^m\subset\cup_{k=1}^{2m}P_k$.
In view of the Finite Union Theorem it is sufficient to show that
$asdim P_k\le n$ for all $k$.
We proceed
by induction on $k$.  It is easy to see that 
$P_{k+1}^A\subset P_k^BA$.
Assuming the inequality $asdim P_k\le n$, we show that $asdim P_{k}^BA\le n$.
We define $Y_r=P_kN^A_r(C)$ where $N^A_r(C)$ denotes an $r$-neighborhood of
$C$ in $A$. First we show that $Y_r\subset N_r(P_k)$. Let $y\in Y_r$, then
$y$ has the form $uz$ where $u\in P^B_k$, $z\in A$ and $ dist(z,C)\le r$, i.e.
$\|z^{-1}c\|\le r$ for some $c\in C$.
Let $c'\bar x_1\dots\bar x_k$ be the normal presentation of $u$,
then $uz=c'x_1x_2\dots x_{k-1}x_kz$ where $x_k\in B\setminus  C$. 
We note that the element $uc$ has the normal presentation
$c'\bar x_1\dots\overline{x_kc}$ and hence $uc\in P_k$. Then
$ dist(y,uc))=\|z^{-1}c\|\le r$, therefore $dist(y,P_k)\le r$, i.e. 
$y\in N_r(P_k)$.
Since the $r$-neighborhood $N_r(P_k)$ is coarsely isomorphic to
the space $P_k$, by the induction assumption we have $asdim N_r(Y_r)\le n$
and hence, $asdim Y_r\le n$.

We consider families $xA$ with $x\in P^B_k$. Let $xA$ and $x'A$ be two different 
cosets.
Since $x$ and $x'$ are different elements with $l(x)=l(x')$, and 
$x^{-1}x'\notin A$, 
the normal presentation of $a^{-1}x^{-1}x'a'$ ends by the 
coset $Ca'$.

Then by Proposition 3 
$
 dist(xA\setminus  Y_r,x'A\setminus  Y_r)=
\|a^{-1}x^{-1}x'a'\|\ge \bar d(Ca',C)= dist(Ca',C)= dist(a',C)>r.
$ 
Note that $P_k^BA$ is the union of these sets $xA$. Since all $xA$ are isometric,
we have a uniform inequality $asdim xA\le n$.
According to Theorem 1 we obtain that $asdim P_k^BA\le n$ and hence 
$asdim P^A_{k+1}\le n$.
Similarly one obtains the inequality $asdim P^B_{k+1}\le n$. The Finite Union
Theorem implies that $asdim P_{k+1}\le n$.
\endprf

{\bf Proof of Theorem 5}\qua
We define a graph $T$ as follows. The vertices of $T$ are the left cosets
$xA$ and $yB$. Two vertices $xA$ and $yB$ are joined by an edge if there is $z$
such that $xA=zA$ and $yB=zB$. 
To check that $T$ is a tree we introduce the weight of a vertex $Y\in T$
given by $w(Y)=\min\{l(y)\mid y\in Y\}$.
Note that for every vertex $e$ with $w(e)>0$ there is a unique neighboring vertex 
$e_-$ with $w(e_-)<w(e)$.
Since we always have $w(zA)\ne w(zB)$, we get
an orientation on $T$ with $w(e_-)<w(e_+)$ for every edge $e$.
The existence this orientation
implies that $T$ does not contain cycles.
Since every vertex of $T$ can be connected with the vertex $A$, the graph $T$
is connected. Thus, $T$ is a tree.
The action of $A\ast_CB$ on $T$ is defined by
left multiplication. We note that the $k$-stabilizer $W_k(A)$ is contained in $(BA)^k$.
Then by Proposition 4 $asdim W_k(A)\le n$. By Theorem 2 $asdim A\ast_CB\le 2n+1$.
\endprf

Let $\{A_i,\|\ \|_i\}$ be a sequence of groups with norms and let 
$C\subset A_i$ be a common subgroup. These norms define a norm $\|\ \|$ on
the amalgamated product $\ast_CA_i$ by taking $\|x\|$ equal the minimum of sums
$\Sigma_{k=1}^{l}\|a_{i_k}\|_{i_k}$ where $x=a_{i_1}\dots a_{i_l}$ and 
$a_{i_k}\in A_{i_k}$.

The following theorem generalizes Theorem 4 and Theorem 5.

\proclaim{Theorem 6}
Let $\{A_i,\|\ \|_i\}$ be a sequence of groups 
satisfying $asdim A_i\le n$ uniformly 
 and let $\|\ \|$ be the norm on a free product
$\ast A_i$ generated by the norms $\|\ \|_i$. Let $C$ be a common subgroup. 
Then $asdim(\ast_C A_i,\|\ \|)\le 2n+1$.
\endproclaim

The proof is omitted since it follows exactly the same scheme.

The following fact will be used in Section 6 in the case of the free 
product.

\proclaim{Proposition 5}
Assume that the groups $A_i$ are supplied with the norms which 
generate the norm on the amalgamated product $\ast_C A_i$.
Let $\psi:\ast_C A_i\to\Gamma$ be a monomorphism to a finitely
generated group such that the restriction $\psi|_{A_i}$ is an isometry 
for every $i$.
Then $\psi$ is a coarsely uniform embedding.
\endproclaim
\prf
Since $\psi$ is a bijection onto the image,
both maps $\psi$ and $\psi^{-1}$ are coarsely proper. We check that
both are coarsely uniform.
First we show that  $\psi$ is 1-Lipschitz. Let $x,y\in\ast_CA_i$
and let $x^{-1}y=a_{i_1}\dots a_{i_n}$ with 
$\|x^{-1}y\|=\Sigma_{k=1}^n\|a_{i_k}\|_{i_k}$. Then
$d_{\Gamma}(\psi(x),\psi(y))\le $\\
$d_{\Gamma}(\psi(x),\psi(xa_{i_1}))+
d_{\Gamma}(\psi(xa_{i_1}),\psi(xa_{i_1}a_{i_2}))+\dots
+d_{\Gamma}(\psi(xa_{i_1}\dots a_{i_{n-1}}),\psi(y))$\\
$=\Sigma_{k=1}^n\|\psi(a_{i_k})\|_{\Gamma}=\Sigma_{k=1}^n\|a_{i_k}\|_{i_k}=
\|x^{-1}y\|= dist(x,y).$\\

Now we show that $\psi^{-1}$ is uniform.
For every $r$ the preimage $\psi^{-1}(B_r(e))$ is finite, since $B_r(e)$ is
finite and $\psi$ is injective. We define 
$\xi(r)=\max\{\|z\|\mid z\in \psi^{-1}(B_r(e))\}$. Let $\bar\xi$ be
strictly monotonic function which tends to infinity and $\bar\xi\ge\xi$.
Let $\rho$ be the inverse function of $\bar\xi$. Then

$ 
d_{\Gamma}(\psi(x),\psi(y))=\|\psi(x^{-1}y)\|_{\Gamma}=
\rho(\bar\xi(\|\psi(x^{-1}y)\|_{\Gamma}))\ge
\rho(\xi(\|\psi(x^{-1}y)\|_{\Gamma}))\ge$\\
$\rho(\|x^{-1}y\|)=\rho(d(x,y))$\\

The last inequality follows from the
inequality $\xi(\|\psi(z)\|)\ge\|z\|$ and the fact 
that $\rho$ is an increasing
function.\endprf

\section{HNN extension}

Let $A$ be a subgroup of a group $G$ and let $\phi: A\to G$ be a monomorphism.
We denote by $G'$ the HNN extension of $G$ by means of $\phi$, i.e. a
group $G'$ generated by $G$ and an element $y$ with the relations 
$yay^{-1}=\phi(a)$ for all $a\in A$. 

\proclaim{Theorem 7} Let $\phi:A\to G$ be a monomorphism of a subgroup
$A$ of a group $G$ with $asdim G\le n$ and let $G'$ be the HNN extension of $G$.
Then $asdim G'\le 2n+1$.
\endproclaim
We recall that a reduced presentation of
an element $x\in G'$ is a word 
$$
g_0y^{\epsilon_1}g_1\dots y^{\epsilon_n}g_n=x,
$$
where $g_i\in G$, $\epsilon_i=\pm 1$, with the property that $g_i\notin A$ whenever
$\epsilon_i=1$ and $\epsilon_{i+1}=-1$ and $g_i\notin \phi(A)$ whenever
$\epsilon_i=-1$ and $\epsilon_{i+1}=1$. The number $n$ is called the length
of the reduced presentation $g_0y^{\epsilon_1}g_1\dots y^{\epsilon_n}g_n$.

The following facts are well-known [\ref{S}]:

A) ({\it uniqueness}) 
Every two reduced presentations of the same element have the same length
and can be obtained from each other by a sequence of the following 
operations:

(1) replacement of $y$ by $\phi(a)ya^{-1}$,\\ 
(2) replacement of $y^{-1}$ by $a^{-1}y\phi(a)$, $a\in A$\\

B) ({\it existence})
Every word of type $g_0y^{\epsilon_1}g_1\dots y^{\epsilon_n}g_n$ can be deformed
to a reduced form by a sequence of the following operations:

(1) replacement of $ygy^{-1}$ by $\phi(g)$  for $g\in A$,
(2) replacement of $y^{-1}\phi(g)y$ by $g$ for $g\in A$,
(3) replacement of $g'\bar g$ by $g=g'\bar g\in G$ if $g',\bar g\in G$.

In particular the uniqueness implies that for any two reduced presentations
$g_0y^{\epsilon_1}g_1\dots y^{\epsilon_n}g_n$ and
$g'_0y^{\epsilon'_1}g_1\dots y^{\epsilon'_n}g'_n$ of the same element
$x\in G'$ we have $(\epsilon_1,\dots,\epsilon_n)=
(\epsilon'_1,\dots,\epsilon'_n)$. 

Let $G$ be a finitely generated group and let $S$ be a finite 
set of generators. We consider the norm on $G'$ defined by the
generating set $S'=S\cup\{y,y^{-1}\}$.

\proclaim{Proposition 6}
Let $g_0y^{\epsilon_1}g_1\dots y^{\epsilon_n}g_n$ be a reduced presentation
of $x\in G'$. Then $\|x\|\ge d(g_n,A)$ if $\epsilon_n=1$ and
$\|x||\ge d(g_n,\phi(A))$ if $\epsilon_n=-1$.
\endproclaim
\prf
We consider here the case when $\epsilon_n=1$. A shortest presentation of $x$
in the alphabet $S'$ gives rise an alternating presentation
$x=r^0_0y^{\epsilon^0_1}r^0_1\dots y^{\epsilon^0_{m_0}}r^0_{m_0}$,
$r_i^0\in G$, $\epsilon^0_i=\pm 1$ with $\|x\|=m_0+\|r^0_0\|+\dots+
\|r^0_{m_0}\|$. We consider
a sequence of presentations of $x$ connecting the above
presentation with a reduced presentation 
$r^1_0y^{\epsilon^1_1}r^1_1\dots y^{\epsilon^1_{m_1}}r^1_{m_1}$
by means of operations (1)-(3) of B).
Then by A) we have that $m_1=n$, $\epsilon_n^1=\epsilon_n=1$ and 
$g_n=\tilde ar_n^1$, $\tilde a\in A$.
Because of the nature of transformations (1)-(3) of B), we can trace out 
to the shortest presentation the letter $y=y^{\epsilon^1_n}$ 
from the reduced word. This means that the 0-th word has the form
$r^0_0y^{\epsilon^0_1}r^0_1\dots y^{\epsilon^0_{l}}r^0_lyw$ where $w$ is an 
alternating word representing $r_{m_k}^k$. Then $\|x\|\ge\|w\|=\|r_{m_k}^k\|=
\|\tilde a^{-1}g_n\|\ge d(g_n,A)$.\endprf

We denote by $l(x)$ the length of a reduced presentation of $x\in G'$.
Let $P_l=\{x\in G\mid l(x)=l\}$.
\proclaim{Proposition 7}
Suppose that $asdim G\le n$, $n>0$. Then $asdim P_l\le n$ for all $l$.
\endproclaim
\prf
We use induction on $l$. We note that $P_0=G$ and
$P_l\subset P_{l-1}yG\cup P_{l-1}y^{-1}G$. We show first that 
$asdim(P_l\cap P_{l-1}yG)\le n$.
Let $r$ be given. We define $Y_r=P_{l-1}yN_r(A)$ where $N_r(A)$ is the 
$r$-neighborhood of $A$ in $G$. We check that 
$Y_r\subset N_{r+1}(P_{l-1})$. Let $z\in Y_r$, then $z=xyg=xyaa^{-1}g=
x\phi(a)ya^{-1}g$ where $x\in P_{l-1}$, $g\in N_r(A)$ and 
$a\in A$ with $\|a^{-1}g\|\le r$. Then $x\phi(a)\in P_{l-1}$ and
$d(x\phi(a),z)=\|ya^{-1}g\|\le\|y\|+\|a^{-1}g\|=r+1$. Since $Y_k$ is coarsely
isomorphic to $P_{l-1}$, by the induction assumption we have $asdim Y_k\le n$.
We consider the family of sets $xyG$ with $x\in P_{l-1}$. If $xyG\ne x'yG$,
then $y^{-1}x^{-1}x'y\notin G$. A reduction in this word can occur only in the
middle. Therefore $x^{-1}x'\notin\phi(A)$. Moreover the reduced presentation
of $y^{-1}x^{-1}x'y$ after these reductions in the middle will be of the form
$y^{-1}r_1\dots r_sy$. Then $d(xyG\setminus  Y_r,x'yG\setminus  Y_r)=
d(xyg,x'yg')=\|g^{-1}y^{-1}x^{-1}x'y'g'\|$. Since
$g^{-1}y^{-1}x^{-1}x'y'g'$ is a reduced presentation, by Proposition 6
$\|g^{-1}y^{-1}x^{-1}x'y'g'\|\ge d(g,A)>r$. So, all the conditions of Theorem 1
are satisfied and, hence $asdim P_{l-1}yG\le n$. Similarly one can show that
$asdim(P_l\cap P_{l-1}y^{-1}G)\le n$. Then the inequality 
$asdim P_l\le n$ follows from the Finite Union Theorem.
\endprf

{\bf Proof of Theorem 7}\qua
We consider a graph $T$ with vertices the left cosets $xG$.
A vertex $xG$ is joined by an edge with a vertex $xgy^{\epsilon}G$, $g\in G$,
$\epsilon=\pm 1$ whenever both $x$ and $xgy^{\epsilon}$ are reduced 
presentations. Since $l(x)=l(xg)$ for all $g\in G$, we can define the
length of a vertex $xG$ of the graph. Thus all edges in $T$ are given
an
orientation and every vertex is connected by a path with the vertex $G$.
Since the length of vertices grows along the orientation, there are no oriented
cycles in $T$. We also note that no vertex can be the end point of two
different edges. All this implies that $T$ is a tree. The group $G'$ acts on $T$
by multiplication from the left. We note that the $r$-stabilizer 
$W_r(G)$ is contained in $P_r$.
Hence by Proposition 7 $asdim W_r(G)\le n$. Then Theorem 2 implies that
$asdim G'\le 2n+1$.\endprf

\rk{Remark}  Both the amalgamated product and the HNN extension are 
the fundamental groups of the simplest graphs of the group [\ref{S}].
We note that theorems of Sections 4-5 can be extended to
the fundamental groups of general graph of groups, since all of them are acting
on the trees with the $R$-stabilizers having an explicit description.

\section{Davis' construction}

We recall that a rightangled Coxeter group W is a group given by the
following presentation:
$$
W=\langle s\in S\mid s^2=1, (ss')^2=1, (s,s')\in E\rangle
$$
where $S$ is a finite set and $E\subset S\times S$.
A barycentric subdivision $N$ of any finite polyhedron defines a rightangled
Coxeter group by the rule: $S=N^{(0)}$ and $E=\{(s,s')\mid (s,s')\in N^{(1)}\}$.
The complex $N$ is called the nerve of $W$ (see [\ref{D}]). We recall that the group
$W$ admits a proper cocompact action on the Davis complex $X$ which is formed as
the union $X=\cup_{w\in W}wC$, where $C=cone(N)$ is called the chamber. Note that the action
of $W$ on the set of centers of the chambers (i.e. cone vertices) is transitive. 
The orbit space of this action is
$C$,  and all isotropy groups are finite. Note that the Davis complex $X$ is contractible.
There is a finite index subgroup $W'$ in $W$ for which the complex
$X/W'$ is a classifying space. 
We denote $\partial C=N$.
Let $X^{\partial}$ denote a subcomplex
$X=\cup_{w\in W}w\partial C\subset X$.
 In [\ref{D}] it was shown that there is a linear 
order on $W$, $e\le w_1\le w_2\le w_3\le \dots$
such that the union $X_{n+1}^{\partial}=\cup_{i=1}^{n+1}w_i\partial C$ is 
obtained  by attaching $w_{n+1}\partial C$ to $X^{\partial}_n$ 
along a contractible subset.
Assume that $N\subset M$ is a subset of an aspherical complex $M$.
We can build the space $X^M$ with an action of the group $W$ on it 
by attaching a copy of $M$ to each $w\partial C$.
Then by induction one can show that every 
complex $X^M_n$ is aspherical and therefore $X^M$  is aspherical.

For every group $\pi$ with $K=K(\pi,1)$ a finite complex, M. Davis considered the
following manifold. Let $M$ be a regular neighborhood of $K\subset{\bf R}^k$ in some 
Euclidean space and let $N$ be a barycentric subdivision of a triangulation of 
the boundary of $M$. Then Davis' manifold is the orbit space $X^M/W'$.
It is aspherical, since $X^M$ is aspherical. 
We refer to the fundamental group  $\Gamma=\pi_1(X^M/W')$ as  
{\it Davis' extension} of the group $\pi$.
By taking a sufficiently large $k$, we may assume that the 
inclusion $N\subset M$ induces an isomorphism
of the fundamental groups. Then in the above notation 
$\Gamma=\pi_1(X^{\partial}/W')$.
\proclaim{Theorem 8}
If $asdim\pi<\infty$, then $asdim\Gamma<\infty$.
\endproclaim
\prf
Since $X^{\partial}$ is path connected, the inclusion $X^{\partial}\subset X$
induces an epimorphism $\phi:\Gamma=\pi_1(X^{\partial}/W')\to\pi_1(X/W')=W'$.
Let $K$ be the kernel. We note that $K=\pi_1(X^M)=\pi_1(X^{\partial})=
\lim_{\to}\{\ast_i\pi_1(w_i\partial C)\}$. It was proven in [\ref{DJ}] that
$asdim W<\infty$. 
The following lemma and Theorem 3 complete the proof.\endprf

\proclaim{Lemma 1}
Assume that $K\subset\Gamma$ is supplied with the induced metric from $\Gamma$.
Then $asdim K\le asdim \pi$.
\endproclaim
\prf We fix a finite generating set $S$ for $\Gamma$.
We consider $A_w=\pi_1(w\partial C)$, $w\in W$ as a subgroup of $K$ defined
by a fixed path $I_w$ joining $x_0$ with $w(x_0)$. Assume that $A_w$ is 
supplied with the norm induced from $\Gamma$. We show that the inequality
$asdim A_w\le asdim\pi$ holds uniformly and by Theorem 4 we obtain that
$asdim(\ast_wA_w,\|\ \|)\le asdim\pi$ for the norm $\|\ \|$ generated by the norms 
on $A_w$. Then we complete the proof applying Proposition 5.

Let $p:X^{\partial}\to\partial C$ be projection onto the orbit space under
the action of $W$. Then $p=q\circ p'$ where $p':X^{\partial}\to X^{\partial}/W'$
is a covering map. We consider the norm on $\pi=\pi_1(\partial C)$ defined by
the generating set
$q_*(S)$. This turns $\pi$ into a metric space of bounded geometry.
Then the homomorphism 
$q_*:\pi_1(X^{\partial})=\Gamma\to\pi_1(\partial C)=\pi$ is 1-Lipschitz map.
The restriction of $q_*$ onto $A_w$ defines an isomorphism acting by 
conjugation with an element generated by the loop $p(I_w)$.
Then according to Proposition 1 we have the inequality $asdim A_w\le asdim\pi$ 
uniformly.
\endprf

%%%%%%%%%%%%%%%%%%%%  End of main body of article
%
\references         %  Print the references here

\Addresses\recd

\bye